\documentclass[a4paper,11pt]{article}
\usepackage{amsmath}
\usepackage{amsfonts}
\usepackage{amssymb}
\usepackage{latexsym}
\usepackage{epsfig}
\usepackage{graphicx}
\usepackage{oldgerm}
\usepackage{theorem}

\def\C{\mathbb{C}}
\def\R{\mathbb{R}}
\def\N{\mathbb{N}}
\def\Z{\mathbb{Z}}
\def\A{\mathbb{A}}
\def\W{\mathbb{W}}
\def\I{\mathbb{I}}
\def\P{{\mathbb P}}
\def\E{{\mathbb E}}
\def\mbK{\mathbb{K}}

\def\x{\mib{x}}
\def\y{\mib{y}}
\def\u{\mib{u}}
\def\v{\mib{v}}
\def\f{\mib{f}}
\def\t{\mib{t}}
\def\X{\mib{X}}
\def\Y{\mib{Y}}
\def\b{\mib{b}}
\def\B{\mib{B}}
\def\n{\mib{n}}

\def\rP{{\rm P}}
\def\rE{{\rm E}}
\def\rC{{\rm C}}

\def\cF{{\cal F}}
\def\cL{{\cal L}}
\def\cH{{\cal H}}
\def\cM{{\cal M}}
\def\cD{{\cal D}}
\def\cK{{\cal K}}

\def\1{{\bf 1}}
\def\0{{\bf 0}}

\def\bE{{\bf E}}
\def\bP{{\bf P}}

\def\mM{\mathfrak{M}}
\def\mS{\mathfrak{S}}
\def\mB{\mathfrak{B}}

\def\supp{{\rm supp}\,}

\def\sp{{\sf p}}
\def\sA{{\sf A}}
\def\sD{{\sf D}}
\def\sS{{\sf S}}

\theorembodyfont{\itshape}

\newtheorem{thm}{Theorem}[section]
\newtheorem{lem}[thm]{Lemma}

\newtheorem{prop}[thm]{Proposition}

\newtheorem{df}[thm]{Definition}
\newtheorem{conj}[thm]{Conjecture}

\newcommand{\mib}[1]{\mbox{\boldmath $#1$}}
\newcommand{\SSC}[1]{\section{#1}\setcounter{equation}{0}}
\newcommand{\qed}{\hbox{\rule[-2pt]{3pt}{6pt}}}

\begin{document}

\title{\bf 
Elliptic Bessel processes and elliptic Dyson models \\
realized as temporally inhomogeneous processes
}
\author{
Makoto Katori \footnote{katori@phys.chuo-u.ac.jp} \\
{\it Department of Physics,
Faculty of Science and Engineering, Chuo University, } \\
{\it Kasuga, Bunkyo-ku, Tokyo 112-8551, Japan}}
\date{22 September 2016}
\pagestyle{plain}
\maketitle
\begin{abstract}
The Bessel process with parameter $D>1$ and the Dyson model
of interacting Brownian motions with coupling constant $\beta >0$
are extended to the processes in which the drift term and the 
interaction terms are given by the logarithmic derivatives of
Jacobi's theta functions.
They are called the elliptic Bessel process, eBES$^{(D)}$,
and the elliptic Dyson model, eDYS$^{(\beta)}$, respectively.
Both are realized on the circumference of a circle 
$[0, 2 \pi r)$ with radius $r >0$
as temporally inhomogeneous processes
defined in a finite time interval $[0, t_*), t_* < \infty$.
Transformations of them to 
Schr\"odinger-type equations with time-dependent potentials
lead us to proving that eBES$^{(D)}$ and eDYS$^{(\beta)}$ can be
constructed as the time-dependent Girsanov transformations
of Brownian motions.
In the special cases where $D=3$ and $\beta=2$,
observables of the processes are defined and 
the processes are represented for them using the Brownian paths
winding round a circle and pinned at time $t_*$. 
We show that eDYS$^{(2)}$
has the determinantal martingale representation for any observable.
Then it is proved that 
eDYS$^{(2)}$ is determinantal for all observables 
for any finite initial configuration
without multiple points.
Determinantal processes are stochastic integrable systems
in the sense that all spatio-temporal correlation functions are
given by determinants controlled by a single continuous
function called the spatio-temporal correlation kernel. 
\end{abstract}


\normalsize

\SSC
{Introduction \label{sec:introduction}}

The transition probability density (tpd) of the standard Brownian motion (BM)
on $\R$ is given by
$
p^{\rm BM}(t, y|x)
=\1(t > 0) e^{-(y-x)^2/2t}/\sqrt{2 \pi t} + \1(t=0)\delta (y-x)$,
$x, y \in \R$,
where $\1(\omega)$ is the indicator function;
$\1(\omega)=1$, if the condition $\omega$ is satisfied, and $\1(\omega)=0$,
otherwise.
If we consider a BM on a circle with radius $r >0$;
$\sS^1(r)=\{x \in \R: x + 2 \pi r=x\}$,
its tpd is obtained by wrapping $p^{\rm BM}$ as
\begin{equation}
p^r_{+1}(t, y|x)=\sum_{w \in \Z} p^{\rm BM}(t, y+2 \pi r w|x),
\quad x, y \in \sS^1(r), \quad t \geq 0,
\label{eqn:pr+1}
\end{equation}
where $w$ denotes the winding number of BM path on the circle. 
Let $i=\sqrt{-1}$. 
Using the Jacobi theta function $\vartheta_3$ defined by
(\ref{eqn:theta}) in Appendix \ref{sec:appendixA},
this is expressed as
\begin{eqnarray}
p^r_{+1}(t, y|x) &=& 
p^{\rm BM}(t, y|x) \vartheta_3 \left( \frac{i (y-x) r}{t} ; \frac{2 \pi i r^2}{t} \right)
\nonumber\\
&=& \frac{1}{2 \pi r} \vartheta_3 \left( \frac{y-x}{2 \pi r}; \frac{it}{2 \pi r^2} \right).
\label{eqn:pr+2}
\end{eqnarray}
Here Jacobi's imaginary transformation 
(\ref{eqn:Jacobi_imaginary}) was performed in the second equality.
The paper \cite{BPY01} is useful to know 
the probability-theoretical interpretations
of formulas for the Jacobi theta functions (and the Riemann zeta function).
See also Sections 13 and 14 of \cite{Bel61}. 
In the present paper, we consider the following `partner' of $p^r_{+1}$,
\begin{equation}
p^r_{-1}(t, y|x) = \sum_{w \in \Z} (-1)^w p^{\rm BM}(t, y+ 2 \pi r w|x),
\quad x, y \in \sS^1(r), \quad t \geq 0,
\label{eqn:pr-1}
\end{equation}
which is expressed as
\begin{eqnarray}
p^r_{-1}(t, y|x) 
&=& p^{\rm BM}(t, y|x) \vartheta_0 \left( \frac{i(y-x)r}{t}; \frac{2 \pi i r^2}{t} \right)
\nonumber\\
&=& \frac{1}{2 \pi r} \vartheta_2 \left( \frac{y-x}{2 \pi r}; \frac{it}{2 \pi r^2} \right).
\label{eqn:pr-2}
\end{eqnarray}
In particular, if we set $y=\pi r$, we have 
\begin{equation}
p^r_{-1}(t, \pi r| x) = \frac{1}{2 \pi r} \vartheta_1 \left( \frac{x}{2 \pi r}; \frac{i t}{2 \pi r^2} \right).
\label{eqn:pr-3}
\end{equation}

In spite of alternating signs in the summation (\ref{eqn:pr-1}),
$p^r_{-1}(t,y|x)$ has appeared in the following problems
in probability theory.

\vskip 0.3cm
\noindent
{\bf Forrester's even-$N$ problem of vicious walkers on a circle} \\

For $N \in \{2,3, \dots\}$, consider the Weyl chamber of type A$_{N-1}$,
$$
\W_N=\{\x=(x_1, x_2, \dots, x_N) \in \R^N : x_1 < x_2 < \cdots < x_N \}.
$$
The tpd of absorbing BM in $\W_N$, in which absorbing walls are
put at the boundaries $\partial \W_N$, is given by \cite{Gra99}
\begin{equation}
q_{\W_N}(t, \y|\x)=\det_{1 \leq j, k \leq N} [p^{\rm BM}(t, y_j|x_k)], \quad 
\x, \y \in \W_N, \quad t \geq 0.
\label{eqn:KM1}
\end{equation}
This determinantal expression is called the Karlin-McGregor formula 
in probability theory \cite{KM59},
and the Lindstr\"om-Gessel-Viennot formula
in enumerative combinatorics \cite{Lin73,GV85}.
Define the Weyl alcove of type A$_{N-1}$ with scale $2 \pi r$ as
\cite{Gra02,Kra07}
\begin{equation}
\A^{2 \pi r}_N=\{\x \in \R^N: x_1 < x_2 < \cdots < x_N < x_1+2 \pi r \},
\label{eqn:alcove}
\end{equation}
and write the tpd of absorbing BM in $\A^{2 \pi r}_N$ as
$q_{\A^{2 \pi r}_N}(t, \y|\x)$, $\x, \y \in \A^{2 \pi r}_N$, $t \geq 0$.
Since the tpd of a single BM on $\sS^1(r)$ is given by (\ref{eqn:pr+1}),
it is expected that $q_{\A^{2 \pi r}_N}(t, \y|\x)$ is given by
$\det_{1 \leq j, k \leq N}[p^r_{+1}(t, y_j|x_k)]$.
It is not the case, however, and if $\pi^{\ell}_N$ represents
the cyclic permutation which shifts indices $\{1,2, \dots, N\}$ by $\ell$, that is, 
$\pi^{\ell}_N(k)=k+\ell \, (\mbox{mod $N$}), 1 \leq k \leq N$, 
and
$\pi^{\ell}_N(\y) \equiv (y_{\pi^{\ell}_N(1)}, \dots, y_{\pi^{\ell}_N(N)})$, 
we find
\begin{equation}
\det_{1 \leq j, k \leq N} [p^r_{+1}(t, y_j|x_k)]
=\sum_{\ell=0}^{N-1} {\rm sgn} (\pi^{\ell}_N)
q_{\A^{2 \pi r}_N}(t, \pi^{\ell}_N(\y)|\x),
\quad \x, \y \in \A^{2 \pi r}_N, \quad t \geq 0.
\label{eqn:detp+1}
\end{equation}
Since ${\rm sgn}(\pi^{\ell}_N) \equiv 1$ for $N$ odd,
$\det_{1 \leq j, k \leq N}[p^r_{+1}(t, y_j|x_k)]$ does not
give a single $q_{\A^{2 \pi r}_N}(t, \y|\x)$, but gives the sum
of tpd's over shifts \cite{For90a},
$$
\sum_{\ell=0}^{N-1} q_{\A^{2 \pi r}_N}(t, \pi^{\ell}_N(\y)|\x)
=\det_{1 \leq j, k \leq N} [p^r_{+1}(t, y_j|x_k)],
\quad \mbox{if $N$ is odd}.
$$
Forrester claimed that, for ${\rm sgn}(\pi^{\ell}_N)=(-1)^{\ell}$
for $N$ even, the formula obtained from (\ref{eqn:detp+1}) is not
useful for $N$ even by alternative signs, and addressed a problem
to find a determinantal formula for 
$\sum_{\ell=0}^{N-1} q_{\A^{2 \pi r}_N}(t, \pi^{\ell}_N(\y)|\x)$ with $N$ even.
This problem was first solved by Fulmek for the non-intersecting lattice-path model
(the vicious walker model) \cite{Ful04} and then by 
Liechty and Wang for the noncolliding BM \cite{LW13} as follows,
$$
\sum_{\ell=0}^{N-1} q_{\A^{2 \pi r}_N}(t, \pi^{\ell}_N(\y)|\x)
=\det_{1 \leq j, k \leq N} [p^r_{-1}(t, y_j|x_k)],
\quad \mbox{if $N$ is even}, 
$$
where $p^r_{-1}$ is given by (\ref{eqn:pr-1}). 
See \cite{NF03,Kat14} for applications of
Forrester-Fulmek formula to prove that
the trigonometric extension of Dyson's BM with $\beta=2$
(the $N$-particle system of noncolliding BMs on $\sS^1(r)$)
is determinantal.

\vskip 0.5cm
\noindent
{\bf Villat's kernel and Komatu-Loewner evolution in an annulus} \\

For $0 < q < 1$, define an annulus
$$
\sA_q=\{z \in \C : q < |z| < 1 \}.
$$
{\it Villat's function} is introduced as
$$
\cK_q(z)= \lim_{N \to \infty} \sum_{n=-N}^N 
\frac{1+q^{2n} z}{1-q^{2n} z}, \quad
z \in \sA_q.
$$
For $z \in \sA_q$ and $\zeta \in \partial \sA_q$, define
{\it Villat's kernel} by
$$
\cK_q(z, \zeta)=\cK_q(z/\zeta).
$$
Suppose that a function $f$ is analytic on $\sA_q$,
continuous on $\overline{\sA_q}$, and $\Re f$ is equal to a
real constant $A$ on 
$\partial \sD_q=\{q e^{i \theta} : 0 \leq \theta < 2 \pi\}=\sS^1(q)$.
Then 
$(1/2 \pi) \int_0^{2 \pi} \Re f(e^{i \theta}) d \theta =A$ and
$f$ can be expressed as
$$
f(z)=\frac{1}{2 \pi} \int_0^{2 \pi} \Re f(e^{i \theta})
\cK_q(z, e^{i \theta}) d \theta + i c,
\quad z \in \sA_q,
$$
with some real constant $c$
(see \cite{FK14} and references therein). 
Fix an annulus $\sA_Q$ with $0 < Q < 1$ and
a Jordan curve $\gamma=\{\gamma(t) : 0 \leq t \leq t_{\gamma} \}$
satisfying $\gamma(0) \in \partial \sD=\sS^1(1)$,
$\gamma(0, t_{\gamma}] \subset \sA_Q$.
Then there exists a strictly increasing function
$\alpha: [0, t_{\gamma}] \mapsto [Q, Q_{\gamma}]$
with $\alpha(t_{\gamma})=Q_{\gamma} < 1$,
which gives the following.
If $\alpha(t)=q \in [Q, Q_{\gamma}]$, then there exists a unique conformal map
$g_q$ from $\sA_Q \setminus \gamma(0, t]$ 
into $\sA_q$ with the normalization condition
$g_q(Q)=q$.
It is proved in \cite{FK14} that
$\alpha$ is continuous, 
$\widetilde{\gamma}(q)=\gamma(\alpha^{-1}(q))$ is well-defined,
$g_q(z), z \in \sA_Q \setminus \gamma(0, t]$ is continuously
differentiable  in $q \in [Q, Q_{\gamma}]$, 
and $g_q(z)$ satisfies the following differential equation
\begin{equation}
\frac{\partial \log g_q(z)}{\partial \log q}
=\cK_q(g_q(z), \lambda(q)) - i \Im \cK_q(q, \lambda(q)),
\quad Q \leq q \leq Q_{\gamma},
\quad q_Q(z)=z,
\label{eqn:KL1}
\end{equation}
where $\lambda(q)=g_q(\widetilde{\gamma}(q))$.
The equation (\ref{eqn:KL1}) is called the
{\it Komatu-Loewner equation} for an annulus.
Its stochastic versions 
were introduced by letting $\lambda(q)$ be 
a time change of BM on
$\partial \sD$ and 
have been extensively studied
in \cite{Zhan04,BF06,BF08,FK14,CFR15}.
Here we claim that, if we set
$$
q=\alpha(t)=e^{-(t_{\gamma}-t)/2}, \quad 0 \leq t \leq t_{\gamma},
$$
then Villat's function is related with $p^1_{-1}$ by
$$
\cK_q(e^{i x}) = 2 i \frac{\partial}{\partial x} \log p^1_{-1}(t_{\gamma}-t, \pi |x),
\quad x \in [0, 2 \pi).
$$

\vskip 0.3cm

Let $0 < t_* < \infty$.
In the present paper, we introduce two families of
stochastic processes.
The first one consists of elliptic extensions of
the $D$-dimensional Bessel process, BES$^{(D)}$. 
It is a one-parameter ($D>1$) 
family of temporally inhomogeneous processes
defined in a finite time-interval $t \in [0, t_*)$.
We consider the stochastic differential equation (SDE) on $\R$,
\begin{equation}
\check{X}(t)=u+ B(t) + \frac{D-1}{2} \int_0^t 
\frac{\partial}{\partial x} 
\log p_{-1}^{r}(t_*-s, \pi r| \check{X}(s)) ds, \quad t \in [0, t_*), 
\label{eqn:SDE_eBES0}
\end{equation}
started at $\check{X}(0)=u \in (0, 2 \pi r)$,
where $B(t), t \geq 0$ is the one-dimensional standard BM
started at the origin.
Here the integrand means
$\partial f(x)/\partial x |_{x=\check{X}(s)}$ with
$f(x)=\log p_{-1}^r(t_*-s, \pi r|x)$. 
Then the process $X(t), t \in [0, t_*)$ is defined on the circumference
$[0, 2 \pi r)$ of $\sS^1(r)$  by
\begin{equation}
X(t)=\check{X}(t) \quad 
(\mbox{mod $2 \pi r$}), \quad t \in [0, t_*).
\label{eqn:eBES2}
\end{equation}
As explicitly shown in Remark 1 given in Sec. \ref{sec:eBES_Bessel},
if we take the double limit $t_* \to \infty$, $r \to \infty$,
then the equation (\ref{eqn:SDE_eBES0}) with (\ref{eqn:eBES2}) is reduced
to the well-known SDE which defines the
$D$-dimensional Bessel process on $[0, \infty)$.
We call this one-parameter family of processes $\{X(t) : t \in [0, t_*)\}$
the {\it elliptic Bessel process with parameter $D > 1$}
and write it
as eBES$^{(D)}, D >1$ for short.

The second one consists of elliptic extensions of
Dyson's Brownian motion model with parameter $\beta>0$ \cite{Dys62}.
It is the one-parameter family ($\beta>0$)
of temporally inhomogeneous 
interacting particle systems
defined in a finite time-interval $[0, t_*)$.
For $N \in \N$, 
express its even-odd parity by
\begin{equation}
\sp(N) = (-1)^N = \left\{ \begin{array}{ll}
+1, \qquad \mbox{if $N$ is even}, \cr
-1, \qquad \mbox{if $N$ is odd}. 
\end{array} \right.
\label{eqn:pN}
\end{equation}
For $\x=(x_1, \dots, x_N) \in \A^{2 \pi r}_N$, we introduce the notation
\begin{equation}
\overline{x} \equiv \sum_{j=1}^N x_j -\kappa_N
\label{eqn:xbar}
\end{equation}
with
\begin{equation}
\kappa_N=\frac{\pi r}{2}(2N-3+\sp(N))
= \left\{ \begin{array}{ll}
\pi r (N-1), \qquad \mbox{if $N$ is even}, \cr
\pi r (N-2), \qquad \mbox{if $N$ is odd}. 
\end{array} \right.
\label{eqn:kappaN}
\end{equation}
We consider the following set of SDEs for an $N$-particle system, 
$\check{\X}(t)=(\check{X}_1(t), \dots, \check{X}_N(t)), t \in [0, t_*)$
on $\R$,
\begin{eqnarray}
\check{X}_j(t) &=& u_j+B_j(t)
+\frac{\beta}{2} \sum_{\substack{1 \leq k \leq N, \cr k \not=j}}
\int_0^t 
\frac{\partial}{\partial x_j} 
\log p^r_{-1}(N(t_*-s), \pi r |\check{X}_j-\check{X}_k) ds
\nonumber\\
&+& \frac{\beta}{2} 
\int_0^t 
\frac{\partial}{\partial x_j} 
\log p^r_{-1}(N(t_*-s), \pi r |\overline{\check{X}}(s)) ds,
\quad 1 \leq j \leq N, \quad t \in [0, t_*),
\label{eqn:eDYS0}
\end{eqnarray}
started at $\u=(u_1, \dots, u_N) \in \A^{2 \pi r}_N$ with $\overline{u} \in (0, 2 \pi r)$, 
where $B_j(t), t \geq 0, 1 \leq j \leq N$ are
independent standard BMs on $\R$. 
We then define the process 
$\X(t)=(X_1(t), \dots, X_N(t)) \in [0, 2 \pi r)^N, t \in [0, t_*)$ by 
\begin{equation}
X_j(t)=\check{X}_j(t) \quad
(\mbox{mod} \, \, 2 \pi r), \quad
1 \leq j \leq N, \quad t \in [0, t_*).
\label{eqn:eDYS3}
\end{equation}
As explicitly shown in Remark 3 given 
at the end of Sec. \ref{sec:eDYS_beta},
if we take the double limit $t_* \to \infty$, $r \to \infty$,
then the system (\ref{eqn:eDYS0}) with (\ref{eqn:eDYS3}) is reduced
to the well-known system of SDEs which defines 
Dyson's Brownian motion model with parameter $\beta$ on $\R$.
We call this one-parameter family of processes $\{ \X(t) : t \in [0, t_*) \}$
the {\it elliptic Dyson model with parameter $\beta > 0$}
and write it as eDYS$^{(\beta)}, \beta >0$ for short.
In the previous paper \cite{Kat15},
we studied eBES$^{(3)}$ and eDYS$^{(2)}$.
Here we treat them as the special cases of the
above two families of processes.

The Kolmogorov equation for eDYS$^{(\beta)}$, $\beta > 0$
is given in the form
$$
\left[ \frac{\partial}{\partial t} - \cL_N^{(\beta)}(t) \right] u(t, \x)=0
$$
with the Kolmogorov operator
\begin{equation}
\cL_N^{(\beta)}(t) = \frac{1}{2} \sum_{j=1}^N \frac{\partial^2}{\partial x_j^2}
+ \frac{\beta}{2} \sum_{j=1}^N \frac{\partial}{\partial x_j}
\frac{\partial W_N}{\partial x_j}(t, \x),
\quad t \in [0, t_*), \quad \x \in \A^{2 \pi r}_N, 
\label{eqn:LN_1}
\end{equation}
associated with a time-dependent potential $W_N(t, \x)$.
(The explicit expression of $W_N(t, \x)$ is given 
by (\ref{eqn:WN_1}) with (\ref{eqn:cW1}).)
If $W_N$ does not depend on time, 
such an $N$-particle diffusion process
is able to be transformed into a Calogero-Moser-Sutherland quantum system
(see Chapter 11 of \cite{For10}).
In the present paper, 
we consider the transformation from 
the temporally inhomogeneous Kolmogorov equation
into a Schr\"odinger-type equation with 
time-dependent Hamiltonian,
\begin{equation}
e^{\beta W_N(t, \x)/2}
\left[ \frac{\partial}{\partial t} -\cL_N^{(\beta)}(t) \right]
e^{-\beta W_N(t, \x)/2}
=\frac{\partial}{\partial t} + (\cH_N^{(\beta)}(t)-E_{N, 0}^{(\beta)}(t)),
\label{eqn:transformationN_1}
\end{equation}
where $\cH_N^{(\beta)}(t)$ has a time-dependent interaction
\begin{equation}
\cH_N^{(\beta)}(t)=-\frac{1}{2} \sum_{j=1}^N \frac{\partial^2}{\partial x_j^2}+V_N^{(\beta)}(t, \x),
\label{eqn:HN_1}
\end{equation}
and $E_{N, 0}^{(\beta)}(t)$ is the ground energy,
which also depends on $t \in [0, t_*)$.
We will show (Lemma \ref{thm:VN}) that 
\begin{eqnarray}
V_N^{(\beta)}(t, \x)
&=& \frac{\beta(\beta-2)}{8} N
\left[
\sum_{1 \leq j < k \leq N}
\left( \zeta \left( x_k-x_j \left| 
2 N \omega_3(t) \right. \right)
-\frac{x_k-x_j}{\pi r} 
\eta_1 \left( 2 \pi r, 2 N \omega_3(t) \right) \right)^2 \right.
\nonumber\\
&& \qquad \qquad 
+ \left( \zeta \left( \overline{x} \left| 
2 \pi r, 2 N \omega_3(t) \right. \right)
-\frac{\overline{x}}{\pi r} 
\eta_1 \left( 2 \pi r, 2 N \omega_3(t) \right) \right)^2 
\nonumber\\
&& \qquad 
\left. -\frac{N-2}{N} \sum_{1 \leq j < k \leq N}
\left( \wp \left( x_k-x_j \left|
2 \pi r, 2 N \omega_3(t) \right. \right) 
+\frac{1}{\pi r} \eta_1 \left( 2 \pi r, 2 N \omega_3(t) \right) \right)
\right],
\nonumber\\
&& \qquad \qquad \qquad \qquad \qquad \qquad 
t \in [0, t_*),  \quad \x \in [0, 2 \pi r)^N,
\label{eqn:VN1b}
\end{eqnarray}
where 
\begin{equation}
2 \omega_3(t) = \frac{i (t_*-t)}{r},
\label{eqn:setting1}
\end{equation}
$\zeta$ and $\wp$ are the Weierstrass functions defined by
(\ref{eqn:zeta1}) and (\ref{eqn:wp1}),
and $\eta_1$ is a special value of $\zeta$ as given by (\ref{eqn:eta1}) 
in Appendix \ref{sec:appendixA}.

Quantum integrable systems have been classified by 
Olshanetsky and Perelomov \cite{OP83}
in connection with Lie algebras.
In their study, interactions were assumed to be time independent
and as the type IV the following interaction 
described using the Weierstrass $\wp$-function was listed, 
\begin{equation}
V_N^{\rm OP_{IV}}(\x)
=g(g-1) \sum_{1 \leq j < k \leq N} \wp(x_k-x_j | 2 \pi r, 2 \omega_3),
\quad \x \in \A^{2 \pi r}_N,
\label{eqn:OP1}
\end{equation}
with a coupling constant $g$ and a constant period $\omega_3$ with $\Im \omega_3 > 0$. 
The many-body quantum system with the interaction (\ref{eqn:OP1})
is called the {\it elliptic Calogero-Moser-Sutherland model}
and has been extensively studied \cite{FV97,Tak00}.
The present potential (\ref{eqn:VN1b}) is time dependent
and different from (\ref{eqn:OP1}).
Although (\ref{eqn:VN1b}) is much complicated, 
appearance of the factor $\beta(\beta-2)$ as the result
of the transformation (\ref{eqn:transformationN_1}) 
implies that the present eDYS$^{(\beta)}$ is a natural extension
of the original Dyson model and its trigonometric/hyperbolic
extensions
(see, for instance, Propositions 11.3.1 and 11.3.2,
and Exercises 11.3.1 in \cite{For10}). 

When $\beta=2$, the interaction (\ref{eqn:VN1b}) vanishes
and the quantum system is allowed to execute a free motion
in a symmetric space $\A_N^{2 \pi r}$.
In the present paper, using the elliptic determinantal equalities of 
Rosengren and Schlosser \cite{RS06}, we will clarify that
this special situation gives the {\it determinantal martingale representation (DMR)} 
\cite{Kat14} to eDYS$^{(2)}$. 
The DMR proves that eDYS$^{(2)}$ is 
a {\it stochastic integrable system} realized as a {\it determinantal process}, 
in the sense that all spatio-temporal correlation functions
are expressed by determinants controlled by
a single function called the spatio-temporal 
correlation kernel \cite{BR05,KT10,Kat16_Springer}.
The present result for eDYS$^{(2)}$ generalizes the previous one 
reported in \cite{Kat15} (see Remark 4 at the end of Sec. \ref{sec:eDYS2}).

The paper is organized as follows.
In Sec. \ref{sec:eBES}, the Kolmogorov equation of eBES$^{(D)}$
is transformed into a one-dimensional Schr\"odinger-type equation
with a potential $V_1^{(D)}$ (Lemma \ref{thm:V1}).
Then we give a probability-theoretical interpretation of this transformation.
That is, we show that eBES$^{(D)}$ is obtained as 
a time-dependent Girsanov transformation \cite{RY99} of a BM (Proposition \ref{thm:Girsanov1}).
When $D=3$, $V_1^{(3)} \equiv 0$
and the transformation can be regarded as a time-dependent extension
of $h$-transformation \cite{RY99}.
We define a collection of observables for eBES$^{(3)}$
(Definition \ref{thm:observable1}) and a signed process
$Y_{-1}(t), t \in [0, t_*)$ in the interval $[0, 2 \pi r)$,
whose finite-dimensional distributions are expressed
using $p_{-1}^r$ given by (\ref{eqn:pr-2})
(Definition \ref{thm:process_Y-1}).
Then, it is proved that expectation of any observable
for eBES$^{(3)}$ is equal to the corresponding expectation
for the process $Y_{-1}$ pinned at the point $\pi r$
at time $t_*$ (Lemma \ref{thm:p-form1}).
In Sec. \ref{sec:eDYS}, the transformation (\ref{eqn:transformationN_1})
from eDYS$^{(\beta)}$ to 
a temporally inhomogeneous version of 
the elliptic Calogero-Moser-Sutherland model
is studied (Lemma \ref{thm:VN} and Proposition \ref{thm:tpd}).
The probability-theoretical realization of this transformation
is the time-dependent Girsanov transformation,
which proves that the probability law of eDYS$^{(\beta)}$ is absolutely
continuous with respect to the Wiener measure for independent BMs
and its Radon-Nikodym derivative gives a martingale (Proposition \ref{thm:GirsanovN}).
The solvability of eDYS$^{(2)}$ is fully discussed in Sec. \ref{sec:eDYS2}.
A collection of observables for eDYS$^{(2)}$ is introduced 
(Definition \ref{thm:observableN})
and $N$-variate signed processes $\Y_{\sp(N)}(t), t \in [0, t_*)$
in $\A^{2 \pi r}_N$, $N \in \{2,3,\dots\}$, are defined by specifying
the finite-dimensional distributions depending on
the parity (\ref{eqn:pN}) of $N$ 
(Definition \ref{thm:process_Ykappa}).
The equivalence is proved between expectation of any observable
for eDYS$^{(2)}$ and that for the process $\Y_{\sp(N)}$
pinned at the $N$-particle configuration $\v=(v_1, \dots, v_N)$ 
with equidistant spacing on $[0, 2 \pi r)$
(Lemma \ref{thm:q-formN1}), where
\begin{equation}
v_j = \frac{\pi r}{2N} (4j-3+\sp(N))
= \left\{ \begin{array}{ll}
\displaystyle{\frac{\pi r}{N}(2j-1)}, \qquad & \mbox{if $N$ is even}, \cr
& \cr
\displaystyle{\frac{2 \pi r}{N}(j-1)}, \qquad & \mbox{if $N$ is odd}, 
\end{array} \right.
\quad
1 \leq j \leq N. 
\label{eqn:v1}
\end{equation}
Here we note that for any $N \in \N$, 
$\overline{v} = \pi r \in (0, 2 \pi r)$.
Then DMR for eDYS$^{(2)}$ is derived for the observables
(Lemma \ref{thm:martingale} and Proposition \ref{thm:DMR}).
Following the definitions of determinantal processes
and spatio-temporal correlation kernels (Definition \ref{thm:determinantal}),
we prove that, for any finite initial configuration without multiple points, 
eDYS$^{(2)}$ is determinantal for all observables 
and we determine the spatio-temporal correlation kernel (Theorem \ref{thm:kernel}). 
Concluding remarks are given in Sec. \ref{sec:concluding_remarks}.

\SSC
{Elliptic Bessel Processes \label{sec:eBES}}

Let $2 \omega_1=2 \pi r$ and using (\ref{eqn:setting1}) define
\begin{equation}
\tau(t) = \frac{\omega_3(t)}{\omega_1}=
\frac{i (t_*-t)}{2 \pi r^2},
\quad t \in [0, t_*).
\label{eqn:setting}
\end{equation}
For a bounded differentiable function $f(t,x) \in \rC_{\rm b}^{1,2}$, we will write
$$
\dot{f}(t,x)=\frac{\partial f}{\partial t}(t,x),
\quad
f^{\prime}(t,x)=\frac{\partial f}{\partial x}(t,x),
\quad
f^{\prime \prime}(t,x)=\frac{\partial^2 f}{\partial x^2}(t,x).
$$
For $N \in \N$ define
\begin{eqnarray}
U_N(t, x)
&=& U_N(t, x; r, t_*) 
\nonumber\\
&=& -\log \vartheta_1 \left(
\frac{x}{2 \pi r}; N \tau(t) \right),
\quad t \in [0, t_*), \quad x \in (0, 2 \pi r).
\label{eqn:cW1}
\end{eqnarray}
The basic properties of this function are given in
Appendix \ref{sec:appendixB}.

\subsection{Elliptic extension of the $D$-dimensional Bessel process
\label{sec:eBES_Bessel}}
Assume $D > 1$.
By (\ref{eqn:pr-3}), we find that 
the SDE (\ref{eqn:SDE_eBES0}) is written using the function $U_1$ as
\begin{equation}
\check{X}(t)=u+ B(t) - \frac{D-1}{2} \int_0^t 
U_1^{\prime}(s, \check{X}(s)) ds, \quad t \in [0, t_*).
\label{eqn:SDE_eBES1}
\end{equation}
As shown by (\ref{eqn:Wdx1}), $U_1^{\prime}(t, x)$ is
a periodic function of $x$ with period $2 \omega_1= 2 \pi r$.
We define 
the elliptic Bessel process with parameter $D > 1$ (eBES$^{(D)}$),
$X(t), t \in [0, t_*)$ on $[0, 2 \pi r)$, by (\ref{eqn:eBES2}).
By the expression (\ref{eqn:Wdx1}), we can see that
\begin{equation}
-U_N^{\prime}(t, x) \sim \left\{
\begin{array}{ll}
\displaystyle{\frac{1}{x}}, 
& \qquad \mbox{as $ x \downarrow 0$}, \cr
\displaystyle{- \frac{1}{2 \pi r-x}},
& \qquad \mbox{as $ x \uparrow 2 \pi r$},
\end{array} \right.
\label{eqn:vicinity}
\end{equation}
independently of $t \in [0, t_*)$ and $N \in \N$.
Therefore, the behavior of $X \in (0, 2 \pi r)$ in the vicinity of 0
(and $2 \pi r$) 
is similar to that of the BES$^{(D)}$ near 0.
This observation implies that if we put
\begin{equation}
T^u_{\{0, 2 \pi r\}} = \inf \{ t > 0: X^u(t) \in \{0, 2 \pi r\} \}
\label{eqn:Tu}
\end{equation}
for eBES$^{(D)}$ started at $u \in (0, 2 \pi r)$, then we can prove that
(see, for instance, \cite{Kat16_Springer}) 
\begin{eqnarray}
\mbox{if $D \geq 2$}, &\mbox{then}& T^u_{\{0, 2 \pi r\}}=\infty,
\forall u \in (0, 2 \pi r), \, \mbox{with probability 1},
\nonumber\\
\mbox{if $1< D < 2$}, &\mbox{then}& T^u_{\{0, 2 \pi r\}} < \infty,
\forall u \in (0, 2 \pi r), \, \mbox{with probability 1}.
\label{eqn:transient_recurrent}
\end{eqnarray}

\vskip 0.3cm
\noindent{\bf Remark 1.} 
Note that the SDE given by (\ref{eqn:SDE_eBES1}) is temporally
inhomogeneous.
The temporally homogeneous system is obtained by
taking the limit $t_* \to \infty$,
\begin{eqnarray}
\check{X}(t) &=& u+ B(t)+ \frac{D-1}{4r} \int_0^t 
\cot \left( \frac{\check{X}(s)}{2r} \right) ds,
\nonumber\\
X(t) &=& \check{X}(t) \quad (\mbox{mod $2 \pi r$}), 
\quad t \geq 0.
\label{eqn:tBES1}
\end{eqnarray}
If we take the further limit $r \to \infty$,
the above becomes
\begin{equation}
X(t)=u+B(t)+\frac{D-1}{2} \int_0^t 
\frac{ds}{X(s)}, \quad t \geq 0.
\label{eqn:BES1}
\end{equation}
The SDE (\ref{eqn:BES1}) defines the $D$-dimensional
Bessel process on $[0, \infty)$, 
and (\ref{eqn:tBES1}) defines its trigonometric extension
defined on $[0, 2\pi r)$.
The process (\ref{eqn:tBES1}) was called a Bessel-like process on $[0, 2 \pi r)$
in Sec. 1.11.1 of \cite{Law05}. 
\vskip 0.3cm

Let $p_1^{(D)}(t, x| t^{\prime}, x^{\prime})$ be the
tpd of eBES$^{(D)}$
from $x^{\prime} \in (0, 2 \pi r)$ at time $t^{\prime} \in [0, t_*)$
to $x \in [0, 2 \pi r]$ at time $t \in (t^{\prime}, t_*)$.
Then the forward Kolmogorov equation (the Fokker-Planck equation)
is defined by
\begin{equation}
\left[ \frac{\partial}{\partial t} - \cL_1^{(D)}(t) \right] p_1^{(D)}(t, x| t^{\prime}, x^{\prime})=0,
\quad t \in [0, t_*), \quad x \in (0, 2 \pi r),
\label{eqn:Kolmogorov1}
\end{equation}
where the time-dependent Kolmogorov operator is defined by
\begin{equation}
\cL_1^{(D)}(t) = \frac{1}{2} \frac{\partial^2}{\partial x^2}
+ \frac{D-1}{2} \frac{\partial}{\partial x} U_1^{\prime}(t, x).
\label{eqn:L1}
\end{equation}
It is well-known that the Kolmogorov operator
can always be transformed into a Hamiltonian operator
when the function $U_1$ is time-independent.
In the present situation, $U_1(t, x)$ is time-dependent
and the Kolmogorov equation including a time-derivative is
transformed into a Schr\"odinger-type equation with 
a time-dependent potential as follows.

\begin{lem}
\label{thm:V1}
For the Kolmogorov operator (\ref{eqn:L1}) of eBES$^{(D)}$, $D > 1$,
the equality
\begin{equation}
e^{(D-1) U_1(t,x)/2}
\left[ \frac{\partial}{\partial t} -\cL_1^{(D)}(t) \right]
e^{-(D-1) U_1(t,x)/2}
= \frac{\partial}{\partial t} + \cH_1^{(D)}(t)
\label{eqn:transformation1}
\end{equation}
is established. 
Here $\cH_1^{(D)}(t)$ is 
the Hamiltonian 
\begin{equation}
\cH_1^{(D)}(t)=-\frac{1}{2} \frac{\partial^2}{\partial x^2}+V_1^{(D)}(t, x).
\label{eqn:H1}
\end{equation}
with a time-dependent potential
\begin{eqnarray}
V_1^{(D)}(t, x) &=& \frac{(D-1)(D-3)}{8}
U_1^{\prime}(t, x)^2
\nonumber\\
&=& \frac{(D-1)(D-3)}{8}
(U_1^{\prime \prime}(t,x)
+ 2 \dot{U}_1(t, x)),
\quad t \in [0, t_*),  \quad x \in (0, 2 \pi r). 
\label{eqn:V1}
\end{eqnarray}
This potential is also written as
\begin{eqnarray}
V_1^{(D)}(t, x) &=& \frac{(D-1)(D-3)}{8}
\left(
\zeta( x | 2 \pi r, 2 \omega_3(t))
- \frac{x}{\pi r} \eta_1( 2 \pi r, 2 \omega_3(t) ) \right)^2,
\nonumber\\
&& \quad t \in [0, t_*), \quad x \in (0, 2 \pi r),
\label{eqn:V1b}
\end{eqnarray}
where $\zeta$ and $\eta_1$ are defined by (\ref{eqn:zeta1}) and
(\ref{eqn:eta1}), respectively, in Appendix \ref{sec:appendixA}.
\end{lem}
\noindent{\it Proof.} \,
For general $f(t,x) \in \rC_{\rm b}^{1,2}$, we find
$$
e^{(D-1)f(t,x)/2} \frac{\partial}{\partial t} e^{-(D-1) f(t,x)/2}
=\frac{\partial}{\partial t} - \frac{D-1}{2} \dot{f}(t, x),
$$
and
\begin{eqnarray}
&& e^{(D-1)f(t,x)/2} 
\left[ \frac{1}{2} \frac{\partial^2}{\partial x^2}
+\frac{D-1}{2} \frac{\partial}{\partial x} f^{\prime}(t, x) \right]
e^{-(D-1) f(t,x)/2}
\nonumber\\
&& \qquad
=\frac{1}{2} \frac{\partial^2}{\partial x^2}
-\frac{(D-1)^2}{8} f^{\prime}(t, x)^2
+\frac{D-1}{4} f^{\prime \prime}(t, x).
\nonumber
\end{eqnarray}
Hence we have (\ref{eqn:transformation1}) and (\ref{eqn:H1}) with
$$
V_1^{(D)}(t,x)=\frac{D-1}{8} \Big\{
(D-1) U_1^{\prime}(t, x)^2
-2 (U_1^{\prime \prime}(t, x)
+2 \dot{U}_1(t, x)) \Big\}.
$$
If we use the equation (\ref{eqn:WNrel}) with $N=1$,
we obtain (\ref{eqn:V1}). 
By (\ref{eqn:Wdx1}), the expression (\ref{eqn:V1b}) is obtained.
The proof is thus completed. \qed
\vskip 0.3cm

Note that, if we omit the term including $\dot{U_1}$
and ignore time dependence in the last expression of (\ref{eqn:V1}),
we have 
$$
V_1^{(D)}(x) \propto U_1^{\prime \prime} = c_1 \wp(x)+ c_2
$$
with some constants $c_1, c_2$ 
and the eigenvalue equation of $\cH_1^{(D)}$
will be identified with Lam\'e's differential equation
(see Sec. 23.4 in \cite{WW27}).
Here we consider the time-dependent potential (\ref{eqn:cW1})
with $N=1$, however, 
and the Hamiltonian
corresponding to eBES$^{(D)}$ is different from 
the Lam\'e operator.

Given $t^{\prime} \in [0, t_*)$, $x^{\prime} \in (0, 2 \pi r)$,
assume that 
$\psi_1^{(D)}(t, x|t^{\prime}, x^{\prime})$ solves the Schr\"odinger-type equation,
\begin{equation}
\left[\frac{\partial}{\partial t}
+\cH_1^{(D)}(t) \right] \psi_1^{(D)}(t, x|t^{\prime}, x^{\prime})=0,
\quad x \in (0, 2 \pi r), \quad t \in (t^{\prime}, t_*),
\label{eqn:psi1_1}
\end{equation}
under the conditions
\begin{eqnarray}
\lim_{t \downarrow t^{\prime}}
\psi_1^{(D)}(t, x|t^{\prime}, x^{\prime}) &=& \delta(x-x'),
\nonumber\\
\psi_1^{(D)}(t, x|t^{\prime}, x^{\prime}) &>& 0.
\label{eqn:psi1_2}
\end{eqnarray}
Then the transformation (\ref{eqn:transformation1}) implies that
\begin{equation}
e^{-(D-1) U_1(t,x)/2} \psi_1^{(D)}(t,x|t^{\prime}, x^{\prime})
=
\left( \vartheta_1 \left(
\frac{x}{2 \pi r}; \tau(t) \right) \right)^{(D-1)/2}
\psi_1^{(D)}(t,x|t^{\prime}, x^{\prime})
\label{eqn:tpd0_1}
\end{equation}
solves the forward Kolmogorov equation of eBES$^{(D)}$.
Since $\vartheta_1(x/2 \pi r; \tau(t))=0$
as $x \downarrow 0$ and $x \uparrow 2 \pi r$, $t \in [0, t_*)$,
and we have assumed $D > 1$,
(\ref{eqn:tpd0_1}) vanishes at $x=0$ and $x=2 \pi r$.
In other words, (\ref{eqn:tpd0_1}) gives
the probability density concentrated on paths of eBES$^{(D)}$,
which stay in the interval $(0, 2 \pi r)$.
Hence we add the condition
\begin{equation}
\lim_{x \downarrow 0} \psi_1^{(D)}(t, x | t^{\prime}, x^{\prime}) =
\lim_{x \uparrow 2 \pi r} \psi_1^{(D)}(t, x| t^{\prime}, x^{\prime})=0
\label{eqn:psi1_3}
\end{equation}
in solving (\ref{eqn:psi1_1}). 
Then the tpd of eBES$^{(D)}$, 
$p_1^{(D)}(t, x| t^{\prime}, x^{\prime})$, $0 < t^{\prime} < t < t_*$,
$x, x^{\prime} \in [0, 2 \pi r)$, satisfies 
\begin{eqnarray}
\left[ \frac{\partial}{\partial t} - \cL_1^{(D)}(t) \right]
p_1^{(D)}(t, x| t^{\prime}, x^{\prime}) &=& 0,
\nonumber\\
\lim_{t \downarrow t^{\prime}}
p_1^{(D)}(t, x| t^{\prime}, x^{\prime}) &=& \delta(x-x^{\prime}),
\label{eqn:tpd1_1}
\end{eqnarray}
is given by
\begin{eqnarray}
p_1^{(D)}(t, x| t^{\prime}, x^{\prime})
&=&\left( \frac{\vartheta_1(x/2 \pi r; \tau(t) )}
{\vartheta_1(x^{\prime}/2 \pi r; \tau(t^{\prime}) )}
\right)^{(D-1)/2}
\psi_1^{(D)}(t, x|t^{\prime}, x^{\prime}),
\nonumber\\
&&
0 < t^{\prime} < t < t_*, \quad
x, x^{\prime} \in (0, 2 \pi r).
\label{eqn:tpd1_2}
\end{eqnarray}

Let $\P^u$ be the probability law of eBES$^{(D)}$
started at $u \in (0, 2 \pi r)$.
Consider the one-dimensional standard BM started at the origin, $b(t), t \geq 0$,
which is independent of $B(t), t \geq 0$ used to give
the SDE (\ref{eqn:SDE_eBES1}) for eBES$^{(D)}$.
We write the probability law of $b(t), t \geq 0$ as $\rP$,
which is called the Wiener measure. 
Put $b^u(t)=u+b(t), t \geq 0$. 
For $x, y \in \R$, $x \wedge y \equiv \min \{x, y\}$. 
Then the following is established.

\begin{prop}
\label{thm:Girsanov1}
{\rm (i)} Assume that $u \in (0, 2 \pi r)$. 
For $D>1$, 
\begin{equation}
M_1^u(t) =
\left( \frac{\vartheta_1(b^u(t)/2 \pi r; \tau(t))}
{\vartheta_1(u/2 \pi r; \tau(0))} \right)^{(D-1)/2}
\exp \left[ - \int_0^t V_1^{(D)}(s, b(s)) ds \right],
\quad t \in [0, t_*), 
\label{eqn:M1}
\end{equation}
is martingale. \\
{\rm (ii)} 
Let 
$$
\sigma^u=\inf \{t >0 : b^u(t) \in \{0, 2 \pi r\}\}.
$$
Then $\P^u$ is absolutely continuous with respect to the Wiener measure
$\rP$ and its Radon-Nikodym derivative with respect to $\rP$
is given by
\begin{equation}
\frac{d \P^u}{d \rP} = M_1^u(t \wedge \sigma^u).
\label{eqn:Radon1}
\end{equation}
\end{prop}
\noindent{\it Proof.} \,
(i) By (\ref{eqn:cW1}) and (\ref{eqn:V1}), (\ref{eqn:M1}) is written as
$$
M_1^u(t)
= \exp \left[ - \frac{D-1}{2} \{ U_1(t, b^u(t))-U_1(0,u)\}
-\frac{(D-1)(D-3)}{8} \int_0^t 
U_1^{\prime}(s, b(s))^2 ds \right].
$$
Assume that $0 \leq t < \sigma^u$.
Since It\^o's formula gives
\begin{eqnarray}
U_1(t, b^u(t))-U_1(0, u)
&=& \int_0^t \dot{U}_1(s, b^u(s)) ds
+ \int_0^t U_1^{\prime}(s, b^u(s)) d b(s)
\nonumber\\
&+& \frac{1}{2} \int_0^t U_1^{\prime \prime}(s, b^u(s)) ds,
\nonumber
\end{eqnarray}
and the equation (\ref{eqn:WNrel}) with $N=1$ holds, 
we have the following expression for $M_1^u(t)$,
$$
M_1^u(t)= \exp \left[ 
-\frac{D-1}{2} \int_0^t U_1^{\prime}(s, b^u(s)) d b(s)
-\frac{1}{2} \left( \frac{D-1}{2} \right)^2
\int_0^t U_1^{\prime}(s, b^u(s))^2 ds \right].
$$
It implies 
\begin{equation}
d M_1^u (t) = - \frac{D-1}{2} U_1^{\prime}(t, b^u(t)) M_1^u(t) 
d b^u(t)
\label{eqn:dM1}
\end{equation}
and hence $M_1^u(t)$ is martingale. \\
(ii) 
Consider the process $\widetilde{b}(t), t \in [0, t_*)$ given by
\begin{equation}
\widetilde{b}(t)=b(t)
+\frac{D-1}{2} \int_0^t U_1^{\prime}(s, b^u(s)) ds,
\quad t \in [0, t_*). 
\label{eqn:Btilde1}
\end{equation}
By It\^o's formula, we see
\begin{eqnarray}
d(\widetilde{b}(t) M_1^u(t))
&=& \widetilde{b}(t) d M_1^u(t)
+M_1^u(t) d \widetilde{b}(t)
+ d \langle \widetilde{b}, M_1^u \rangle_t
\nonumber\\
&=& \widetilde{b}(t) d M_1^u(t)+M_1^u(t) d b(t),
\nonumber
\end{eqnarray}
since (\ref{eqn:dM1}) and (\ref{eqn:Btilde1}) give
$$
d \langle \widetilde{b}, M_1^u \rangle_t
=-\frac{D-1}{2} U_1^{\prime}(t, b^u(t)) M_1^u(t) dt,
$$
and
$$
M_1^u(t) d \widetilde{b}(t)
= M_1^u(t) d b(t)+ M_1^u(t) \frac{D-1}{2} U_1^{\prime}(t, b^u(t)) dt.
$$
Then, with respect to $\rP$,
$\widetilde{b}(t) M_1^u(t)$ is martingale.
Assume that $\P^u$ is the probability measure such that 
its Radon-Nikodym derivative
with respect to $\rP$ is given by (\ref{eqn:Radon1}). 
This is concentrated on paths with $\sigma^u > t$
due to the factor $(\vartheta_1(b^u(t)/2 \pi t; \tau(t)) )^{(D-1)/2}, D > 1$
as proved in (i) above.
Hence, with respect to $\P^u$, $\widetilde{b}(t), t \geq 0$ 
can be regarded as a standard BM. In other words, in $\P^u$, 
(\ref{eqn:Btilde1}) should be read as
$$
b^u(t)=u+\widetilde{b}(t)- \frac{D-1}{2}
\int_0^t U_1^{\prime}(s, b^u(s)) ds,
\quad t \in [0, t_*),
$$
with a standard BM, $\widetilde{b}(t)$.
It is equivalent with the SDE (\ref{eqn:SDE_eBES1}), which has defined
eBES$^{(D)}$.
Thus the proof is completed. \qed
\vskip 0.3cm

The above is a time-dependent extension 
of Girsanov's transformation 
(see, for example, Chapter VIII of \cite{RY99}) 
for the present process.

\subsection{Observables of eBES$^{(3)}$ and process $Y_{-1}$
\label{sec:eBES3}}
In this section we consider eBES$^{(D)}$ with
$$
D=3, 
$$
in which Lemma \ref{thm:V1} gives $V_1^{(3)}(t, x) \equiv 0$.
For eBES$^{(3)}$, $X(t), t \in [0, t_*)$,
the natural filtration is considered,
$\cF_X(t)=\sigma(X(s), s \in [0, t]), t \in [0, t_*)$, 
which satisfies the usual conditions. 
Let $u \in (0, 2 \pi r)$ and $T \in [0, t_*)$.
The expectations with respect to $\P^u$ and $\rP^u$
are written as $\E^u$ and $\rE^u$, respectively.
Proposition \ref{thm:Girsanov1} implies that, for any
$\cF_X(T)$-measurable bounded function $F$,
\begin{equation}
\E^u [ F(X(\cdot))]
=\rE^{u} \left[ F(b(\cdot)) \1(\sigma > T)
\frac{\vartheta_1(b(T)/2 \pi r; \tau(T))}{\vartheta_1(u/2 \pi r; \tau(0))} \right],
\label{eqn:Exp1}
\end{equation}
where $\sigma=\inf\{t >0: b(t) \in \{0, 2 \pi r\} \}$.

For an $\cF_X(T)$-measurable function $F$, 
it will be sufficient to consider the case that
$F$ is given as
\begin{equation}
F(X(\cdot))=\prod_{m=1}^M g_m(X(t_m))
\label{eqn:F_g}
\end{equation}
for an arbitrary $M \in \N$ and for
an arbitrary increasing series of times,
$0 \leq t_1 < \dots < t_M \leq T < t_*$ 
with bounded measurable functions $g_m, 1 \leq m \leq M$.

\begin{df}
\label{thm:observable1}
For $T \in [0, t_*)$, 
define the $\cF_X(T)$-observables
of eBES$^{(3)}$ 
as a set of $\cF_X(T)$-measurable functions
given in the form (\ref{eqn:F_g}) satisfying the following conditions. \\
{\rm (1)} 
$\{g_m(x) \}_{m=1}^M$ are periodic functions, 
\begin{equation}
g_m(x+2 \pi r n)=g_m(x), \quad \forall n \in \Z, \quad 1 \leq m \leq M.
\label{eqn:g_period}
\end{equation}
{\rm (2)} 
$\{g_m(x) \}_{m=1}^M$ are symmetric for reflections 
at the boundaries $\{0, 2 \pi r\}$ of the interval $(0, 2 \pi r)$.
That is, 
\begin{equation}
g_m(-x)=g_m(x), \quad
g_m(2\pi r+x)=g_m(2 \pi r -x), \quad
x \in (0, 2 \pi r), \quad 1 \leq m \leq M.
\label{eqn:g_sym}
\end{equation}
\end{df}

By the reflection principle of BM, if $F$ is an $\cF_X(T)$-observable, 
we can omit the indicator function $\1(\sigma > T)$ in (\ref{eqn:Exp1}),
\begin{equation}
\E^u [ F(X(\cdot))]
=\rE^{u} \left[ F(b(\cdot)) 
\frac{\vartheta_1(b(T)/2 \pi r; \tau(T))}{\vartheta_1(u/2 \pi r; \tau(0))} \right],
\label{eqn:Exp2}
\end{equation}
since $\vartheta_1(x/2 \pi r; \tau(t)), t \in [0, t_*)$
is anti-symmetric at the boundaries $\{0, 2 \pi r\}$
of the interval $(0, 2\pi r)$ in the sense that
\begin{equation}
\vartheta_1 \left(- \frac{x}{2\pi r} ; \tau(t) \right) 
=-\vartheta_1 \left( \frac{x}{2 \pi r} ; \tau(t) \right),
\quad
\vartheta_1 \left( \frac{2\pi r+x}{2\pi r}; \tau(t) \right)
=-\vartheta_1 \left( \frac{2 \pi r-x}{2 \pi r}; \tau(t) \right),
\label{eqn:theta_sym}
\end{equation}
$x \in (0, 2 \pi r), t \in [0, t_*)$.

Now we notice the relation (\ref{eqn:pr-3}) between
$\vartheta_1$ and $p_{-1}^r$, the `wrapped' transition
probability density of BM with alternating signs (\ref{eqn:pr-1}).
We see
\begin{eqnarray}
\frac{\vartheta_1(x/2\pi r; \tau(t))}{\vartheta_1(u/2 \pi r; \tau(0))}
&=& \frac{p_{-1}^r(t_*-t, \pi r|x)}{p_{-1}^r(t_*, \pi r|u)}
\nonumber\\
&=& \sum_{w \in \Z} \mu_1^u(w)
\frac{p^{\rm BM}(t_*-t, \pi r+2 \pi r w|x)}{p^{\rm BM}(t_*, \pi r+2 \pi r w|u)}
\nonumber
\end{eqnarray}
with
$$
\mu_1^u(w)=\frac{(-1)^w p^{\rm BM}(t_*, \pi r + 2 \pi r w|u)}
{\sum_{n \in \Z} (-1)^n p^{\rm BM}(t_*, \pi r+ 2 \pi r n|u)}.
$$

\begin{df}
\label{thm:process_Y-1}
The process $Y_{-1}(t), t \in [0, t_*)$ is defined in the 
interval $[0, 2 \pi r)$ by the following finite-dimensional distributions.
Assume that the initial state is given by $u \in [0, 2 \pi r)$.
For an arbitrary $M \in \N$, 
an arbitrary strictly increasing sequence of times
$0 \leq t_1 < \cdots < t_M < t_*$,
and arbitrary Borel sets $A_m \in \mB([0, 2 \pi r))$, $1 \leq m \leq M$, 
\begin{eqnarray}
&& \bP^u[Y_{-1}(t_1) \in A_1, Y_{-1}(t_2) \in A_2, \dots,
Y_{-1}(t_M) \in A_M]
\nonumber\\
&& \quad
= \prod_{m=1}^M \int_{A_m } dx^{(m)} \,
\prod_{m=1}^M p_{-1}^r(t_m-t_{m-1}, x^{(m)}|x^{(m-1)}), 
\nonumber
\end{eqnarray}
where $t_0=0, x^{(0)}=u$.
The expectation is denoted by $\bE^u$.
\end{df}

\begin{lem}
\label{thm:p-form1}
Let $T \in [0, t_*)$. For any $\cF_X(T)$-observable $F$,
\begin{equation}
\E^u[F(X(\cdot))]
= \bE^u \left[ F(Y_{-1}(\cdot))
\frac{p_{-1}^r(t_*-T, \pi r | Y_{-1}(T))}{p_{-1}^r(t_*, \pi r|u)} \right], 
\quad u \in (0, 2 \pi r).
\label{eqn:p-form1}
\end{equation}
\end{lem}
\noindent{\it Proof.} \,
If the observable $F$ is given by (\ref{eqn:F_g}), (\ref{eqn:Exp2}) is
explicitly written as
\begin{eqnarray}
\E^u[F(X(\cdot))]
&=& \sum_{w \in \Z} \mu_1^u(w)
\prod_{m=1}^M \int_{-\infty}^{\infty} dx^{(m)} 
\prod_{m=1}^M \Big( g_m(x^{(m)}) p^{\rm BM}(t_m-t_{m-1}, x^{(m)}|x^{(m-1)}) \Big)
\nonumber\\
&& \qquad \qquad
\times \frac{p^{\rm BM}(t_*-t_M, \pi r + 2 \pi r w|x^{(M)})}
{p^{\rm BM}(t_*, \pi r + 2 \pi r w|u)}. 
\label{eqn:Exp3}
\end{eqnarray}
It is a summation over $w \in \Z$ of the Brownian bridges
of time duration $t_*$, each of which
is started at $u \in (0, 2 \pi r)$ and ended at
$\pi r+ 2 \pi r w \in \R$.
Note that $\mu_1^u(w)$ is a signed measure on $\Z$ with the property
\begin{equation}
\mu_1^u(w+n)=(-1)^n \mu_1^u(w), \quad
w, n \in \Z.
\label{eqn:mu_alt}
\end{equation}
Due to the summation with this alternative-signed measure,
contributions of paths which do not satisfy the condition
$\sigma > T \geq t_M$ are completely canceled,
and the equivalence between (\ref{eqn:Exp1}) and (\ref{eqn:Exp2})
is guaranteed for any observable $F$.

In the multiple integral (\ref{eqn:Exp3}), we rewrite the integrals as
\begin{eqnarray}
\int_{-\infty}^{\infty} dx^{(m)} \, f(x^{(m)})
&=& \sum_{\ell_m \in \Z} 
\int_{2 \pi r \ell_m}^{2 \pi r(\ell_m+1)} dx^{(m)} \, f(x^{(m)})
\nonumber\\
&=& \sum_{\ell_m \in \Z} \int_0^{2 \pi r} d y^{(m)} \,
f(y^{(m)}+2 \pi r \ell_m), \quad 1 \leq m \leq N, 
\nonumber
\end{eqnarray}
where $f$ represents the integrand.
Note that $g_m$'s in (\ref{eqn:Exp3}) are assumed to be 
periodic as (\ref{eqn:g_period}).
Moreover, note the facts that 
\begin{eqnarray}
&&p^{\rm BM}(t_m-t_{m-1}, y^{(m)}+2 \pi r \ell_m | y^{(m-1)}+2 \pi r \ell_{m-1})
\nonumber\\
&& \qquad = p^{\rm BM}(t_m-t_{m-1}, y^{(m)}+2 \pi r (\ell_m-\ell_{m-1}) | y^{(m-1)}),
\quad 1 \leq m \leq M-1
\nonumber
\end{eqnarray}
with $\ell_0=0$,
\begin{eqnarray}
&& p^{\rm BM}(t_M-t_{M-1}, y^{(M)}+2 \pi r \ell_M | y^{(M-1)}+2 \pi r \ell_{M-1})
p^{\rm BM}(t_*-t_M, \pi r + 2 \pi r w| y^{(M)}+2 \pi r \ell_M)
\nonumber\\
&& =
p^{\rm BM}(t_M-t_{M-1}, y^{(M)}+2 \pi r (\ell_M-\ell_{M-1}) | y^{(M-1)})
p^{\rm BM}(t_*-t_M, \pi r + 2 \pi r (w-\ell_M) | y^{(M)}),
\nonumber
\end{eqnarray}
and 
\begin{eqnarray}
&& \sum_{w \in \Z} \mu_1^u(w) p^{\rm BM}
(t_*-t_M, \pi r + 2 \pi r(w-\ell_M)|y^{(M)})
\nonumber\\
&& \qquad \qquad = \sum_{n \in \Z} \mu_1^u(n+\ell_M)
p^{\rm BM}(t_*-t_M, \pi r+ 2 \pi r n|y^{(M)})
\nonumber\\
&& \qquad \qquad
=(-1)^{\ell_M} \sum_{n \in \Z} \mu_1^u(n)
p^{\rm BM}(t_*-t_M, \pi r+ 2 \pi r n|y^{(M)}), 
\nonumber
\end{eqnarray}
where (\ref{eqn:mu_alt}) was used.
Since $(-1)^{\ell_M}=(-1)^{\ell_1} \prod_{m=2}^M (-1)^{\ell_m-\ell_{m-1}}$,
the following expression is derived through the
definition (\ref{eqn:pr-1}) of $p_{-1}^r$,
\begin{equation}
\E^u[F(X(\cdot))]
=\prod_{m=1}^M \int_0^{2 \pi r} d y^{(m)} \,
\prod_{m=1}^M \Big(
g_m(y^{(m)}) p_{-1}^r(t_m-t_{m-1}, y^{(m)}|y^{(m-1)}) \Big)
\frac{p_{-1}^r(t_*-t_M, \pi r|y^{(M)})}{p_{-1}^r(t_*, \pi r|u)}.
\nonumber
\end{equation}
Thus the proof is completed. \qed

Equation (\ref{eqn:p-form1}) shows that
eBES$^{(3)}$ is realized as the process $Y_{-1}$ pinned
at the point $\pi r$ at time $t_*$.
See, for instance, Sec. IV.4 of \cite{BS02}
for the basic properties of pinned Brownian motions
(Brownian bridges).

\SSC{Elliptic Dyson Models \label{sec:eDYS}}
\subsection{Elliptic extensions of Dyson models with parameter $\beta$
\label{sec:eDYS_beta}}
Let $0 < t_* < \infty, 0 < r < \infty, N \in \{2,3, \dots\}$.
Using (\ref{eqn:cW1}) define
\begin{equation}
W_N(t, \x) = \sum_{1 \leq j < k \leq N} U_N(t, x_k-x_j)
+ U_N(t, \overline{x}), 
\quad t \in [0, t_*), \quad
\x=(x_1, \dots, x_N) \in \A^{2 \pi r}_N.
\label{eqn:WN_1}
\end{equation}
Assume that $\beta >0$. 
The system of SDEs (\ref{eqn:eDYS0}) for $N$-particle system
is then simply written as
\begin{equation}
\check{X}_j(t)=u_j+B_j(t)- \frac{\beta}{2} \int_0^t 
\frac{\partial W_N}{\partial x_j} (s, \check{\X}(s)) ds,
\quad 1 \leq j \leq N, \quad t \in [0, t_*),
\label{eqn:eDYS1}
\end{equation}
where $\u=(u_1, \dots, u_N) \in \A^{2 \pi r}_N$ with $\overline{u} \in (0, 2 \pi r)$, and 
$B_j(t), t \geq 0, 1 \leq j \leq N$ are
independent standard BMs on $\R$.
We define the elliptic Dyson model with parameter $\beta >0$ (eDYS$^{(\beta)}$), 
$\X(t)=(X_1(t), \dots, X_N(t)) \in [0, 2 \pi r)^N, t \in [0, t_*)$ by 
(\ref{eqn:eDYS3}).

Since $U_N^{\prime}(t, x)$ is an odd function of $x$,
the summation of (\ref{eqn:eDYS1}) over $j=1,\dots, N$ gives
\begin{equation}
\overline{\check{X}}(t)
=\overline{u}+\sqrt{N} B(t)
- \frac{N \beta}{2} \int_0^t U_N^{\prime}(s, \overline{\check{X}}(s)) ds
\quad t \in [0, t_*),
\label{eqn:eDYS2}
\end{equation}
where $\overline{\check{X}}(t)=\sum_{j=1}^N \check{X}_j(t) -\kappa_N$,
$\overline{u}=\sum_{j=1}^N u_j -\kappa_N$ with (\ref{eqn:kappaN}),
and $B(t), t \geq 0$ is a standard BM on $\R$
which is independent of $B_j(t), t \geq 0, 1 \leq j \leq N$ in (\ref{eqn:eDYS1}).
If we perform the time change $t \mapsto \widetilde{t}=Nt$ and put
$\widetilde{X}(\, \widetilde{t} \,)=\overline{\check{X}}(t)$, then we have
\begin{equation}
\widetilde{X}(\, \widetilde{t} \,)
= \overline{u}+B(\, \widetilde{t} \,)
-\frac{\beta}{2} \int_0^{\widetilde{t}} 
\widetilde{U}_1^{\prime}(s, \widetilde{X}(s)) ds,
\quad \widetilde{t} \in [0, N t_*),
\label{eqn:eDYS2b}
\end{equation}
where 
$$
\widetilde{U}_1^{\prime}(t, x)
=-\frac{\partial}{\partial x} \log \vartheta_1 \left(
\frac{x}{2 \pi r}; \frac{i(Nt_*-t)}{2 \pi r^2} \right),
\quad t \in [0, Nt_*), \quad x \in (0, 2 \pi r).
$$
That is, $\overline{\check{X}}(t)$ of eDYS$^{(\beta)}$ 
is a time change of eBES$^{(D)}$, $\check{X}(t)$, with $D=\beta+1$
defined in $[0, N t_*)$.
(If we set $N=1$ in (\ref{eqn:kappaN}) and (\ref{eqn:v1}), we have
$\kappa_1=-\pi r$, $v_1=0$, and $\overline{v}=v_1-\kappa_1=\pi r$. 
Then the $N=1$ limit of $\overline{\check{X}}(t)$ of eDYS$^{(\beta)}$
will correspond to $\check{X}(t)+\pi r$, where
$\check{X}(t)$ denotes eBES$^{(D)}$ with $D=\beta+1$.)

By the definition (\ref{eqn:xbar}) with (\ref{eqn:kappaN}), 
$$
\overline{x} \in (0, 2 \pi r)
\Longleftrightarrow 
\pi r \left( 1- \frac{3-\sp(N)}{2N} \right)
< \frac{1}{N} \sum_{j=1}^N x_j
< \pi r \left( 1 + \frac{\sp(N)+1}{2N} \right).
$$
Then we consider the subspace of the Weyl alcove (\ref{eqn:alcove}) 
by adding a restriction on the location of
`center of mass' of the system, $\sum_{j=1}^N x_j/N$, 
\begin{eqnarray}
\check{\A}_N^{2 \pi r} 
&\equiv& \A_N^{2 \pi r} \cap 
\{ \x \in \R^N : \overline{x} \in (0, 2 \pi r) \}
\nonumber\\
&=& \A_N^{2 \pi r} \cap
\left\{ \x \in \R^N : 
\pi r \left( 1- \frac{3-\sp(N)}{2N} \right)
< \frac{1}{N} \sum_{j=1}^N x_j
< \pi r \left( 1 + \frac{\sp(N)+1}{2N} \right) \right\},
\label{eqn:alcove_check}
\end{eqnarray}
with $\partial \check{\A}_N^{2 \pi r}=\partial \A_N^{2 \pi r} \cup \{ \x \in \R^N : \overline{x} \in \{0, 2 \pi r\}\}$
and $\overline{\check{\A}_N^{2 \pi r}} = \check{\A}_N^{2 \pi r} \cup \partial \check{\A}_N^{2 \pi r}$.
Put
\begin{equation}
T^{\u}_{\partial \check{\A}_N^{2 \pi r}}
= \inf \{ t > 0: \X^{\u}(t) \in \partial \check{\A}_N^{2 \pi r} \}
\label{eqn:Tu2}
\end{equation}
for eDYS$^{(\beta)}$ started at $\u \in \check{\A}_N^{2 \pi r}$.
Then by (\ref{eqn:vicinity}), the following 
can be proved by the argument given by \cite{RS93,GM13},
\begin{eqnarray}
\mbox{if $\beta \geq 1$}, &\mbox{then}& T^{\u}_{\partial \check{\A}_N^{2 \pi r}}=\infty,
\forall \u \in \check{\A}_N^{2 \pi r}, \, \mbox{with probability 1},
\nonumber\\
\mbox{if $0< \beta < 1$}, &\mbox{then}& T^{\u}_{\partial \check{\A}_N^{2 \pi r}} < \infty,
\forall \u \in \check{\A}_N^{2 \pi r}, \, \mbox{with probability 1}.
\label{eqn:phase_transition}
\end{eqnarray}

\vskip 0.3cm
\noindent{\bf Remark 2.} 
By the relations between Jacobi's theta functions (\ref{eqn:theta_relations})
and quasi-periodicity (\ref{eqn:quasi_periodic}),  we see that
\begin{eqnarray}
U_N^{\prime}(t, \overline{x}) &=&
-\frac{1}{2 \pi r} \frac{\vartheta_1^{\prime}(\sum_{j=1}^N x_j/2 \pi r -\kappa_N/2 \pi r; N \tau(t))}
{\vartheta_1(\sum_{j=1}^N x_j/2 \pi r -\kappa_N/2 \pi r; N \tau(t))}
\nonumber\\
&=&
-\frac{1}{2 \pi r} \frac{\vartheta_1^{\prime}(\sum_{j=1}^N x_j/2 \pi r +1/2; N \tau(t))}
{\vartheta_1(\sum_{j=1}^N x_j/2 \pi r +1/2; N \tau(t))}
\nonumber\\
&=&
-\frac{1}{2 \pi r} \frac{\vartheta_2^{\prime}(\sum_{j=1}^N x_j/2 \pi r; N \tau(t))}
{\vartheta_2(\sum_{j=1}^N x_j/2 \pi r; N \tau(t))}.
\nonumber
\end{eqnarray}
Then the system of SDEs (\ref{eqn:eDYS1}) is expressed as follows without using $\kappa_N$,
\begin{eqnarray}
\check{X}_j(t) &=& u_j+ B_j(t)
+\frac{\beta}{4 \pi r} \sum_{\substack{1 \leq k \leq N, \cr k \not=j}}
\int_0^t 
\frac{\vartheta_1^{\prime}((\check{X}_j(s)-\check{X}_k(s))/2 \pi r; N \tau(s))}
{\vartheta_1((\check{X}_j(s)-\check{X}_k(s))/2 \pi r; N \tau(s))} ds
\nonumber\\
&& \quad + \frac{\beta}{4 \pi r} \int_0^t 
\frac{\vartheta_2^{\prime}(\sum_{\ell=1}^N \check{X}(s)/2 \pi r; N \tau(s))}
{\vartheta_2(\sum_{\ell=1}^N \check{X}(s)/2 \pi r; N \tau(s))} ds,
\quad 1 \leq j \leq N, \quad t \in [0, t_*).
\label{eqn:eDYS1b}
\end{eqnarray}
\vskip 0.3cm
\vskip 0.3cm
\noindent{\bf Remark 3.} 
The system of SDEs (\ref{eqn:eDYS1b}) is temporally inhomogeneous.
The temporally homogeneous system is obtained by taking
the limit $t_* \to \infty$ as
\begin{eqnarray}
\check{X}_j(t) &=& u_j + B_j(t)+ \frac{\beta}{4r}
\sum_{\substack{1 \leq k \leq N, \cr k \not=j}}
\int_0^t 
\cot \left( \frac{\check{X}_j(s)-\check{X}_k(s)}{2r} \right) ds
- \frac{\beta}{4r} 
\int_0^t
\tan \left( \frac{1}{2r} \sum_{\ell=1}^N \check{X}_{\ell}(s) \right) ds,
\nonumber\\
X_j(t) &=& \check{X}_j(t)
\quad (\mbox{mod} \, \, 2 \pi r), \quad
1 \leq j \leq N, \quad t \geq 0,
\label{eqn:tDYS1}
\end{eqnarray}
where we have used (\ref{eqn:theta_asym}) in Appendix \ref{sec:appendixA}.
This system is a 
trigonometric version of the Dyson model with parameter $\beta >0$.
The SDEs have terms including $\sum_{\ell=1}^N \check{X}_{\ell}(s)$ and 
this system is different from the trigonometric Dyson model
studied in \cite{NF03,Kat14}.
If we take the further limit $r \to \infty$, the above system becomes
Dyson's Brownian motion model with parameter $\beta$ defined on $\R$ \cite{Dys62}, 
\begin{equation}
X_j(t) = u_j+B_j(t) + \frac{\beta}{2} 
\sum_{\substack{1 \leq k \leq N, \cr k \not=j}} \int_0^t 
\frac{ds}{X_j(s)-X_k(s)},
\quad 1 \leq j \leq N, \quad t \geq 0.
\label{eqn:DYS1}
\end{equation}
\vskip 0.3cm

\subsection{Temporally inhomogeneous version
of elliptic Calogero-Moser-Sutherland model
\label{sec:CSmodel}}

Here we perform the transformation from 
the Kolmogorov equation
with the time-dependent potential $W_N(t, \x)$
into the Schr\"odinger-type equation
with the time-dependent Hamiltonian $\cH_N^{(\beta)}(t)$
announced in Sec. \ref{sec:introduction}. 

\begin{lem}
\label{thm:VN}
Assume $N \in \{2,3, \dots\}$. 
For the Kolmogorov operator (\ref{eqn:LN_1}) of the
temporally inhomogeneous eDYS$^{(\beta)}, \beta>0$,
the equality
\begin{equation}
e^{\beta W_N(t, \x)/2}
\left[ \frac{\partial}{\partial t} -\cL_N^{(\beta)}(t) \right]
e^{-\beta W_N(t, \x)/2}
=\frac{\partial}{\partial t} + (\cH_N^{(\beta)}(t)-E_{N, 0}^{(\beta)}(t)),
\label{eqn:transformationN_1b}
\end{equation}
is established. 
Here $\cH_N^{(\beta)}(t)$ is the Hamiltonian
\begin{equation}
\cH_N^{(\beta)}(t)=-\frac{1}{2} \sum_{j=1}^N \frac{\partial^2}{\partial x_j^2}+V_N^{(\beta)}(t, \x)
\label{eqn:HN_1b}
\end{equation}
with a time-dependent interaction
\begin{eqnarray}
V_N^{(\beta)}(t, \x) 
&=&
\frac{\beta(\beta-2)}{8} N
\left[
\sum_{1 \leq j < k \leq N}  
U_N^{\prime}(t, x_k-x_j)^2
+U_N^{\prime}(t, \overline{x})^2
\right.
\nonumber\\
&& \qquad 
\left. 
- \frac{N-2}{N}
\sum_{1 \leq j < k \leq N} U_N^{\prime \prime}(t, x_k-x_j)
\right], \quad t \in [0, t_*),  \x \in [0, 2 \pi r)^N,
\label{eqn:VN1}
\end{eqnarray}
where $U_N(t, x)$ is defined by (\ref{eqn:cW1}), 
and $E_{N, 0}^{(\beta)}(t)$ is
the time-dependent ground energy
\begin{equation}
E_{N, 0}^{(\beta)}(t) = 
- \frac{\beta^2}{16} N(N-1)(N-2)
\frac{1}{\pi r}
\eta_1 ( 2 \pi r, 2 N \omega_3(t) ). 
\quad t \in [0, t_*). 
\label{eqn:E0_1}
\end{equation}
The above potential $V_N^{(\beta)}(t, \x)$ is also expressed as 
\begin{eqnarray}
V_N^{(\beta)}(t, \x)
&=& \frac{\beta(\beta-2)}{8} N
\left[
\sum_{1 \leq j < k \leq N}
\left( \zeta \left( x_k-x_j \left| 
2 N \omega_3(t) \right. \right)
-\frac{x_k-x_j}{\pi r} 
\eta_1 \left( 2 \pi r, 2 N \omega_3(t) \right) \right)^2 \right.
\nonumber\\
&& \qquad \qquad 
+ \left( \zeta \left( \overline{x} \left| 
2 \pi r, 2 N \omega_3(t) \right. \right)
-\frac{\overline{x}}{\pi r} 
\eta_1 \left( 2 \pi r, 2 N \omega_3(t) \right) \right)^2 
\nonumber\\
&& \qquad 
\left. -\frac{N-2}{N} \sum_{1 \leq j < k \leq N}
\left( \wp \left( x_k-x_j \left|
2 \pi r, 2 N \omega_3(t) \right. \right) 
+\frac{1}{\pi r} \eta_1 \left( 2 \pi r, 2 N \omega_3(t) \right) \right)
\right],
\nonumber\\
&& \qquad \qquad \qquad \qquad \qquad \qquad 
t \in [0, t_*),  \quad \x \in [0, 2 \pi r)^N.
\label{eqn:VN1c}
\end{eqnarray}
\end{lem}
\noindent{\it Proof.} \,
For general $W_N(t, \x)$, we find
$$
e^{\beta W_N(t, \x)/2} \frac{\partial}{\partial t}
e^{-\beta W_N(t, \x)/2}
=\frac{\partial}{\partial t}-\frac{\beta}{2} \frac{\partial W_N}{\partial t}(t, \x)
$$
and (see, for example, Exercises 11.1.1 in \cite{For10})
$$
e^{\beta W_N(t, \x)/2} \cL_N^{(\beta)}
e^{-\beta W_N(t, \x)/2}
=\frac{1}{2} \sum_{j=1}^N \frac{\partial^2}{\partial x_j^2}
-\frac{\beta^2}{8} \sum_{j=1}^N 
\left( \frac{\partial W_N}{\partial x_j}(t, \x) \right)^2
+\frac{\beta}{4} \sum_{j=1}^N
\frac{\partial^2 W_N}{\partial x_j^2}(t, \x).
$$
Then we have (\ref{eqn:transformationN_1b}) and (\ref{eqn:HN_1b}) with
\begin{eqnarray}
&& V_N^{(\beta)}(t, \x)-E_{N, 0}^{(\beta)}(t)
\nonumber\\
&& 
\quad
= \frac{\beta}{2} \left[
- \frac{\partial W_N}{\partial t}(t, \x)
+ \frac{\beta}{4} \sum_{j=1}^N \left(
\frac{\partial W_N}{\partial x_j}(t, \x) \right)^2
-\frac{1}{2} \sum_{j=1}^N \frac{\partial^2 W_N}{\partial x_j^2} (t, \x) \right]
\nonumber\\
&& \quad
= f_1(t, \overline{x})+f_2(t, \x)
\label{eqn:VN_c1}
\end{eqnarray}
with
\begin{eqnarray}
f_1(t, \overline{x})
&=& \frac{\beta}{2} \left[
- \dot{U}_N(t, \overline{x})
+\frac{\beta}{4} \sum_{j=1}^N \left(
\frac{\partial U_N}{\partial x_j}(t, \overline{x}) \right)^2
-\frac{1}{2} \sum_{j=1}^N \frac{\partial^2 U_N}{\partial x_j^2}(t, \overline{x}) 
\right],
\nonumber\\
f_2(t, \x)
&=& \frac{\beta}{2} \left[ 
- \sum_{1 \leq k < \ell \leq N} \dot{U}_N(t, x_{\ell}-x_k)
+\frac{\beta}{4} \sum_{j=1}^N \left(
\frac{\partial}{\partial x_j} \sum_{1 \leq k < \ell \leq N} 
U_N(t, x_{\ell}-x_k) \right)^2
\right.
\nonumber\\
&& \qquad \left. 
-\frac{1}{2} \sum_{j=1}^N \frac{\partial^2}{\partial x_j^2}
\sum_{1 \leq k < \ell \leq N} U_N(t, x_{\ell}-x_k) \right],
\nonumber
\end{eqnarray}
where the definition (\ref{eqn:WN_1}) of $W_N(t, \x)$ was used.
Since $\overline{x}=\sum_{j=1}^N x_j-\kappa_N$,
$\partial U_N(t, \overline{x})/\partial x_j=U_N^{\prime}(t, \overline{x})$,
$\partial^2 U_N(t, \overline{x})/\partial x_j^2=U_N^{\prime \prime}(t, \overline{x})$
for any $1 \leq j \leq N$.
Thus, if we use (\ref{eqn:WNrel}), we obtain
\begin{equation}
f_1(t, \overline{x})=\frac{\beta(\beta-2)}{8} N U_N^{\prime}(t, \overline{x})^2.
\label{eqn:f1}
\end{equation}
If we use the symmetries (\ref{eqn:BB}), we can show that
\begin{eqnarray}
&& \frac{\partial}{\partial x_j}
\sum_{1 \leq k < \ell \leq N} U_N(t, x_{\ell}-x_k)
= - \sum_{k: k \not=j} U_N^{\prime}(t, x_k-x_j),
\nonumber\\
&& \frac{\partial^2}{\partial x_j^2}
\sum_{1 \leq k < \ell \leq N} U_N(t, x_{\ell}-x_k)
=\sum_{k: k \not=j} U_N^{\prime \prime}(t, x_k-x_j),
\quad 1 \leq j \leq N,
\nonumber
\end{eqnarray}
and 
\begin{eqnarray}
&&
\sum_{j=1}^N \frac{\partial}{\partial x_j}
\left( \sum_{1 \leq k < \ell \leq N}
U_N(t, x_{\ell}-x_k) \right)^2
\nonumber\\
&& \quad
= 2 \sum_{1 \leq j < k \leq N} U_N^{\prime}(t, x_k-x_j)^2
+ \sum_{\substack{ 1 \leq j, k, \ell \leq N, \cr
j \not=k \not= \ell \not= j}}
U_N^{\prime}(t, x_j-x_k) U_N^{\prime}(t, x_j-x_{\ell}).
\nonumber
\end{eqnarray}
Thus we have
\begin{eqnarray}
f_2(t, \x) &=& 
\frac{\beta}{4} \sum_{1 \leq j < k \leq N}
\Big\{ (N-2) U_N^{\prime \prime}(t, x_k-x_j)
-(N-\beta) U_N^{\prime}(t, x_k-x_j) \Big\}
\nonumber\\
&& +\frac{\beta^2}{8}
\sum_{\substack{ 1 \leq j, k, \ell \leq N, \cr
j \not=k \not= \ell \not= j}}
U_N^{\prime}(t, x_j-x_k) U_N^{\prime}(t, x_j-x_{\ell}).
\label{eqn:f2}
\end{eqnarray}
Using the expression (\ref{eqn:Wdx1}) for $U_N^{\prime}(t, x)$,
we obtain the equality
\begin{eqnarray}
&& \sum_{\substack{1 \leq j, k, \ell \leq N, \cr j \not= k \not= \ell \not=j}}
U_N^{\prime} (t, x_j-x_k)
U_N^{\prime}(t, x_j-x_{\ell})
\nonumber\\
&& = (N-2) \sum_{1 \leq j < k \leq N} 
\zeta(x_j-x_k)^2 - (N-2) \sum_{1 \leq j < k \leq N} \wp(x_j-x_k)
\nonumber\\
&& \quad - \frac{2}{\pi r} \eta_1
\sum_{\substack{1 \leq j, k, \ell \leq N, \cr j \not= k \not= \ell \not= j}}
\zeta(x_j-x_k)(x_j-x_{\ell})
+\left( \frac{\eta_1}{\pi r}  \right)^2
\sum_{\substack{1 \leq j, k, \ell \leq N, \cr j \not= k \not= \ell \not= j}}
(x_j-x_k)(x_j-x_{\ell}),
\nonumber
\end{eqnarray}
where we have omitted writing dependence
of $2 \omega_1=2 \pi r$ and $2 N \omega_3(t)=iN(t_*-t)/r$ in
$\zeta, \wp$ and $\eta_1$ just for simplicity of expressions. 
It is easy to verify that
\begin{eqnarray}
&&
\sum_{\substack{1 \leq j, k, \ell \leq N, \cr j \not= k \not= \ell \not=j}}
\zeta(x_j-x_k)(x_j-x_{\ell})
-(N-2) \sum_{1 \leq j < k \leq N}
\zeta(x_k-x_j)(x_k-x_j)
= 0,
\nonumber\\
&& 
\sum_{\substack{1 \leq j, k, \ell \leq N, \cr j \not= k \not= \ell \not=j}}
(x_j-x_k)(x_j-x_{\ell})
-(N-2) \sum_{1 \leq j < k \leq N}
(x_k-x_j)^2
= 0, 
\nonumber
\end{eqnarray}
using the fact that $\zeta(z)$ is odd.
Therefore, we obtain the equalities
\begin{eqnarray}
&& 
\sum_{\substack{1 \leq j, k, \ell \leq N, \cr j \not= k \not= \ell \not=j}}
U_N^{\prime} (t, x_j-x_k)
U_N^{\prime}(t, x_j-x_{\ell})
\nonumber\\
&& =
(N-2) \left[
\sum_{1 \leq j < k \leq N}
\left( - \zeta(x_k-x_j)+\frac{x_k-x_j}{\pi r} \eta_1 \right)^2
-\sum_{1 \leq j < k \leq N}
\left( \wp(x_k-x_j) + \frac{1}{\pi r} \eta_1 \right)^2 \right]
\nonumber\\
&& \quad
+\frac{1}{2} N (N-1)(N-2) \frac{1}{\pi r} \eta_1
\nonumber\\
&& = 
(N-2) \left[
\sum_{1 \leq j < k \leq N} U_N^{\prime}(t, x_k-x_j)^2
-\sum_{1 \leq j < k \leq N} U_N^{\prime \prime}(t, x_k-x_j) \right]
+\frac{1}{2} N (N-1) (N-2) \frac{1}{\pi r} \eta_1,
\nonumber
\end{eqnarray}
where (\ref{eqn:Wdx1}) and (\ref{eqn:Wd2x1}) were used.
Hence (\ref{eqn:f2}) is written as
\begin{eqnarray}
f_2(t, \x) &=&
\frac{\beta (\beta-2)}{8} N
\left[ \sum_{1 \leq j < k \leq N}
U_N^{\prime}(t, x_k-x_j)^2
-\frac{N-2}{N} \sum_{1 \leq j < k \leq N}
U_N^{\prime \prime}(t, x_k-x_j) \right]
\nonumber\\
&&
+ \frac{\beta^2}{16} N (N-1)(N-2) \frac{1}{\pi r} \eta.
\label{eqn:f2b}
\end{eqnarray}
Inserting (\ref{eqn:f1}) and (\ref{eqn:f2b}) in
(\ref{eqn:VN_c1}), 
(\ref{eqn:transformationN_1b}) is proved
with (\ref{eqn:HN_1b})-(\ref{eqn:E0_1}).
The equalities (\ref{eqn:Wdx1}) and (\ref{eqn:Wd2x1}) verify
the equivalence between (\ref{eqn:VN1}) and (\ref{eqn:VN1c}). 
Thus the proof is completed.
\qed
\vskip 0.3cm
For $t^{\prime} \in [0, t_*)$, $\x^{\prime} \in \check{\A}^{2 \pi r}_N$,
define $\psi_N^{(\beta)}(\cdot, \cdot|t^{\prime}, \x^{\prime})$ 
as the unique solution of the Schr\"odinger-type equation,
\begin{equation}
\left[\frac{\partial}{\partial t} +\cH_N^{(\beta)}
\right] \psi_N^{(\beta)}(t, \x|t^{\prime}, \x^{\prime})=0,
\quad \x \in \check{\A}^{2 \pi r}_N, \quad t \in (t^{\prime}, t_*)
\label{eqn:psiN_1}
\end{equation}
with the Hamiltonian (\ref{eqn:HN_1b}), 
under the conditions
\begin{eqnarray}
\lim_{t \downarrow t^{\prime}}
\psi_N^{(\beta)}(t, \x|t^{\prime}, \x^{\prime}) &=& \delta(\x-\x')  \equiv \prod_{j=1}^N \delta(x_j-x^{\prime}_j),
\nonumber\\
\psi_N^{(\beta)}(t, \x|t^{\prime}, \x^{\prime}) &>& 0, \quad \x \in \A^{2 \pi r}_N,
\nonumber\\
\lim_{\x \to \partial \check{\A}^{2 \pi r}_N} 
\psi_N^{(\beta)}(t, \x|t^{\prime}, \x^{\prime}) &=& 0.
\label{eqn:conditionsA}
\end{eqnarray}

\subsection{Transition probability density of eDYS$^{(\beta)}$
\label{sec:tpd}}

The tpd of eDYS$^{(\beta)}$
from $\x^{\prime} \in \check{\A}^{2 \pi r}_N$
at time $t^{\prime} \in [0, t_*)$
to $\x \in \check{\A}^{2 \pi r}_N$ at time $t \in (t^{\prime}, t_*)$,
is given by the unique solution of
\begin{equation}
\left[ \frac{\partial}{\partial t} - \cL_N^{(\beta)}(t) \right]
p_N^{(\beta)}(t, \x| t^{\prime}, \x^{\prime}) = 0,
\label{eqn:tpdN_1}
\end{equation}
under the condition
$$
\lim_{t \downarrow t^{\prime}}
p_N^{(\beta)}(t, \x| t^{\prime}, \x^{\prime}) = \delta(\x-\x^{\prime}),
$$
where $\cL_N^{(\beta)}(t)$ is given by (\ref{eqn:LN_1}).

For $N \in \N, t \geq 0, \x, \x^{\prime} \in \R^N$, define
\begin{equation}
q^r_N(t, \x|\x^{\prime})
=\det_{1 \leq j, k \leq N}
[ p^r_{\sp(N)} (t, x_j|x^{\prime}_k) ], 
\label{eqn:qrN}
\end{equation}
where 
\begin{equation}
p^r_{\sp(N)}(t, x|x^{\prime})
= \sum_{w \in \Z} \sp(N)^w p^{\rm BM}(t, x+2 \pi r w|x^{\prime}),
\quad t \geq 0, 
\label{eqn:pkappa}
\end{equation}
with the parity (\ref{eqn:pN}) of $N$.

\begin{prop}
\label{thm:tpd}
\begin{description}
\item{\rm (i)} \,
Let
\begin{equation}
h_N^r(t, \x)
=\eta(N \tau(t))^{-\beta(N-1)(N-2)/4}
\vartheta_1 \left( \frac{\overline{x}}{2 \pi r};
N \tau(t) \right)
\prod_{1 \leq j < k \leq N}
\vartheta_1 \left( \frac{x_k-x_j}{2 \pi r};
N \tau(t) \right),
\label{eqn:hN_1}
\end{equation}
$t \in [0, t_*), \x \in \R^N$,
where $\eta$ is the Dedekind modular function 
defined by (\ref{eqn:Dedekind1}). Then
\begin{equation}
p_N^{(\beta)}(t, \x| t^{\prime}, \x^{\prime})=
\left( 
\frac{h_N^r(t, \x)}{h_N^r(t^{\prime}, \x^{\prime})}
\right)^{\beta/2}
\psi_N^{(\beta)}(t, \x|t^{\prime}, \x^{\prime}),
\label{eqn:tpdN_2}
\end{equation}
$0 < t^{\prime} < t < t_*, 
\x^{\prime} \in \check{\A}^{2 \pi r}_N,  \x \in \overline{\check{\A}_N^{2 \pi r}}$,
where $\psi_N^{(\beta)}(t, \x|t^{\prime}, \x^{\prime})$
is the unique solution of (\ref{eqn:psiN_1}) under the conditions (\ref{eqn:conditionsA}).
\item{\rm (ii)} The function (\ref{eqn:hN_1}) is also written as 
\begin{equation}
h_N^r(t, \x)
=\eta(N \tau(t))^{-(\beta-2)(N-1)(N-2)/4}
\left( \frac{2 \pi r}{\sqrt{N}} \right)^N
q^r_N(t_*-t, \v|\x), 
\quad t \in [0, t_*),
\label{eqn:hN_2}
\end{equation}
where $\v=(v_1, \dots, v_N)$ is 
the $N$-particle configuration (\ref{eqn:v1}) with equidistant spacing
on $[0, 2 \pi r)$.
\end{description}
\end{prop}
\noindent{\it Proof.} \,
(i) 
The transformation (\ref{eqn:transformationN_1b}) 
and (\ref{eqn:psiN_1}) implies
that 
$$
\left[ \frac{\partial}{\partial t}
-\cL_N^{(\beta)}(t) + E_{N, 0}^{(\beta)}(t) \right]
e^{-\beta W_N(t, \x)/2} 
\psi_N^{(\beta)}(t, \x|t^{\prime}, \x^{\prime})
=0.
$$
Then 
$$
C(t) e^{-\beta W_N(t, \x)/2} 
\psi_N^{(\beta)}(t, \x|t^{\prime}, \x^{\prime})
$$
with
$$
C(t)=\exp \left[
\int_0^t E_{N, 0}^{(\beta)}(s) ds \right]
$$
solves the forward Kolmogorov equation (\ref{eqn:tpdN_1}).
The factor $C(t)$ is the unique solution of the differential equation
\begin{equation}
\frac{d}{dt} \log C(t)=E_N^{(\beta)}(t), \quad t \in [0, t_*)
\label{eqn:logC}
\end{equation}
with the initial condition $C(0)=1$.
By (\ref{eqn:Dedekind2}),  (\ref{eqn:E0_1}) given in 
Lemma \ref{thm:VN} is written as
\begin{eqnarray}
E_N^{(\beta)}(t)
&=& \frac{\beta^2}{16} N(N-1)(N-2) \frac{i}{\pi r^2}
\left. \frac{d}{d \tau} \log \eta(\tau) \right|_{\tau=N \tau(t)}
\nonumber\\
&=& - \frac{\beta^2}{8}(N-1)(N-2) \frac{d}{dt} \log \eta(N\tau(t)),
\nonumber
\end{eqnarray}
where (\ref{eqn:setting1}) and (\ref{eqn:setting}) were used.
That is, we have obtained
$$
E_N^{(\beta)}(t)=\frac{d}{dt} \log \eta(N \tau(t))^{-\beta^2(N-1)(N-2)/8}.
$$
Comparing this with (\ref{eqn:logC}), we can conclude that
\begin{equation}
C(t)=
\exp \left[
\int_0^t E_{N, 0}^{(\beta)}(s) ds \right]
= \left( \frac{\eta(N \tau(t))}{\eta(N \tau(0))} 
\right)^{-\beta^2 (N-1)(N-2)/8},
\quad t \in [0, t_*).
\label{eqn:Dedekind3}
\end{equation}
Since $W_N(t, \x)$ is defined by (\ref{eqn:WN_1}) with
(\ref{eqn:cW1}), we see
$$
e^{-\beta W_N(t, \x)/2}
=\left[ \vartheta_1 \left( \frac{\overline{x}}{2 \pi r};
N \tau(t) \right)
\prod_{1 \leq j < k \leq N}
\vartheta_1 \left( \frac{x_k-x_j}{2 \pi r};
N \tau(t) \right) \right]^{\beta/2}.
$$
Hence (\ref{eqn:tpdN_2}) with (\ref{eqn:hN_1}) is proved. \\
(ii) The following equalities were proved as Proposition 5.6.3
in \cite{For10}. 
Let $\alpha \in \C$. Assume that $\Im \tau > 0$. 
For $N$ odd
\begin{eqnarray}
&& \det_{1 \leq j, k \leq N} 
\left[ \vartheta_3 \left(x_j+\alpha- \frac{k}{N}; \tau \right) \right]
= N^{N/2} \eta(N \tau)^{-(N-1)(N-2)/2}
\nonumber\\
&& \qquad \qquad \times
\vartheta_3 \left( \sum_{j=1}^N (x_j+\alpha); N \tau \right)
\prod_{1 \leq j < k \leq N} \vartheta_1(x_k-x_j; N \tau),
\label{eqn:Forrester1}
\end{eqnarray}
while for $N$ even
\begin{eqnarray}
&& \det_{1 \leq j, k \leq N} 
\left[ \vartheta_1 \left(x_j+\alpha-\frac{k}{N}; \tau \right) \right]
= N^{N/2} \eta(N \tau)^{-(N-1)(N-2)/2}
\nonumber\\
&& \qquad \qquad \times
\vartheta_0 \left( \sum_{j=1}^N (x_j+\alpha); N \tau \right)
\prod_{1 \leq j < k \leq N} \vartheta_1(x_k-x_j; N \tau).
\label{eqn:Forrester2}
\end{eqnarray}
The equalities (\ref{eqn:Forrester1}) and (\ref{eqn:Forrester2})
were given as Lemma 2.6 in our previous paper \cite{Kat15},
but there are errors and
$\vartheta_{\mu}(\sum_{j=1}^N (x_j+\alpha)+N \tau/2; 2 N \tau)$,
$\mu=0, 3$ should be replaced by
$\vartheta_{\mu}(\sum_{j=1}^N (x_j+\alpha); N \tau)$,
as shown in (\ref{eqn:Forrester1}) and (\ref{eqn:Forrester2}). 

For $N$ odd, we put $\alpha=1/N+\tau/2$ in (\ref{eqn:Forrester1}),
and for $N$ even, $\alpha=1/2N+1/2+\tau/2$ in (\ref{eqn:Forrester2}).
Let $\tau=\tau(t)$. Then, 
for (\ref{eqn:qrN}) with (\ref{eqn:pkappa}), 
we obtain the equality
$$
q^r_N(t_*-t, \v|\x)
=\left( \frac{\sqrt{N}}{2 \pi r} \right)^N
\eta(N\tau(t))^{-(N-1)(N-2)/2}
\vartheta_1 \left( \frac{\overline{x}}{2 \pi r};
N \tau(t) \right)
\prod_{1 \leq j < k \leq N}
\vartheta_1 \left( \frac{x_k-x_j}{2 \pi r};
N \tau(t) \right),
$$
where (\ref{eqn:pr+2}), (\ref{eqn:pr-3}), (\ref{eqn:theta_relations})
and (\ref{eqn:quasi_periodic}) 
have been used.
Thus (\ref{eqn:hN_2}) is proved. 
\qed
\vskip 0.3cm

The product formula (\ref{eqn:hN_1}) shows that 
if $\x \in \check{\A}^{2 \pi r}_N$,
then $h_N^r(t, \x)>0$.
On the other hand, the determinantal formula (\ref{eqn:hN_2}) 
shows that, for $p^r_{\sp(N)}(t_*-t, v_j|x), 1 \leq j \leq N$
are not identically zero as functions of $x$,
$h_N^r(t, \x)$ becomes zero only when
$\x$ arrives at $\partial \check{\A}^{2 \pi r}_N$ from the inside of $\check{\A}_N^{2 \pi r}$.
Moreover, if we write $\pi(\x)=(x_{\pi(1)}, \dots, x_{\pi(N)})$
for a permutation $\pi \in \mS_N$ of $N$ indices of 
$\x=(x_1, \dots, x_N)$, the determinantal expression (\ref{eqn:hN_2})
implies
\begin{equation}
h_N^r(t, \pi(\x))=
{\rm sgn}(\pi) h_N^r(t, \x), \quad
t \in [0, t_*), \quad \pi \in \mS_N.
\label{eqn:hN_3}
\end{equation}

\subsection{Time-dependent Girsanov's transformation
\label{sec:Girsanov}}
Let $\P^{\u}$ be the probability law of eDYS$^{(\beta)}$
started at $\u \in \check{\A}^{2 \pi r}_N$.
Consider the $N$-dimensional standard BM started at the origin, 
$\b(t)=(b_1(t), \dots, b_N(t)), t \geq 0$,
which is independent of $\B(t), t \geq 0$ used to give
the SDEs (\ref{eqn:eDYS1}) for eDYS$^{(\beta)}$.
We write the Wiener measure of $\b(t), t \geq 0$ as $\rP$.
Put $\b^{\u}(t)=\u+\b(t), t \geq 0$.
Then as a multivariate extension of Proposition \ref{thm:Girsanov1}, 
the following is proved.

\begin{prop}
\label{thm:GirsanovN}
{\rm (i)} Assume that $\u \in \check{\A}^{2 \pi r}_N$.
For $\beta > 0$, 
\begin{equation}
M_N^{\u}(t) =
\left( \frac{h_N^r(t, \b^{\u}(t))}{h_N^r(0, \u)} \right)^{\beta/2}
\exp \left[ - \int_0^t V_N^{(\beta)}(s, \b^{\u}(s)) ds \right],
\quad t \in [0, t_*),
\label{eqn:MN}
\end{equation}
is martingale. \\
{\rm (ii)}
Let 
\begin{equation}
\sigma^{\u}=\inf \{t >0 : \b^{\u}(t) \in \partial \check{\A}^{2 \pi r}_N \}.
\label{eqn:sN}
\end{equation}
Then $\P^{\u}$ is absolutely continuous with respect to the Wiener measure
$\rP$ and its Radon-Nikodym derivative with respect to $\rP$
is given by
\begin{equation}
\frac{d \P^{\u}}{d \rP} = M_N^{\u}(t \wedge \sigma^{\u}).
\label{eqn:RadonN}
\end{equation}
\end{prop}
\noindent{\it Proof.} \,
(i) We will prove that $M_N^{\u}(t)$ given by (\ref{eqn:MN}) is equal to
\begin{equation}
\widetilde{M}_N^{\u}(t)
= \exp \left[
-\frac{\beta}{2} \sum_{j=1}^N \int_0^t 
\frac{\partial W_N}{\partial x_j}(s, \b^{\u}(s)) db_j(s)
-\frac{1}{2} \left( \frac{\beta}{2} \right)^2
\sum_{j=1}^N \int_0^t 
\left( \frac{\partial W_N}{\partial x_j}(s, \b^{\u}(s) \right)^2 ds \right],
\label{eqn:tildeMN1}
\end{equation}
$t \in [0, t_*)$. 
It\^o's formula gives
\begin{eqnarray}
&& W_N(t, \b^{\u}(t))-W_N(0, \u)
\nonumber\\
&& =
\int_0^t 
\frac{\partial W_N}{\partial s}(s, \b^{\u}(s)) ds
+ \sum_{j=1}^N \int_0^t 
\frac{\partial W_N}{\partial x_j}(s, \b^{\u}(s)) d b_j^{u_j}(s) 
+\frac{1}{2} \sum_{j=1}^N \int_0^t 
\frac{\partial^2 W_N}{\partial x_j^2} (s, \b^{\u}(s)) ds.
\nonumber
\end{eqnarray}
Then, the logarithm of (\ref{eqn:tildeMN1}) is written as
\begin{eqnarray}
\log \widetilde{M}_N^{\u}(t)
&=& -\frac{\beta}{2} \Big\{
W_N(t, \b^{\u}(t))-W_N(0, \u) \Big\}
\nonumber\\
&& + \frac{\beta}{2} \int_0^t 
\frac{\partial W_N}{\partial s}(s, \b^{\u}(s)) ds
+\frac{\beta}{4} \sum_{j=1}^N \int_0^t 
\frac{\partial^2 W_N}{\partial x_j^2} (s, \b^{\u}(s)) ds
\nonumber\\
&&
-\frac{1}{2} \left( \frac{\beta}{2} \right)^2
\sum_{j=1}^N \int_0^t 
\left( \frac{\partial W_N}{\partial x_j}(s, \b^{\u}(s)) \right)^2 ds.
\nonumber
\end{eqnarray}
In the proof of Lemma \ref{thm:VN}, we have shown that
\begin{eqnarray}
&& - \frac{\beta}{2} \frac{\partial W_N}{\partial t}(t, \x)
- \left\{ - \frac{\beta^2}{8}
\sum_{j=1}^N \left( \frac{\partial W_N}{\partial x_j}(t, \x) \right)^2
+\frac{\beta}{4} \sum_{j=1}^N \frac{\partial^2 W_N}{\partial x_j^2}(t, \x) \right\}
\nonumber\\
&& \quad
= e^{\beta W_N(t, \x)/2} 
\left[ \frac{\partial}{\partial t}-\cL_N^{(\beta)} \right] e^{-\beta W_N(t, \x)/2}
-\left( \frac{\partial}{\partial t}
-\frac{1}{2} \sum_{j=1}^N \frac{\partial^2}{\partial x_j^2} \right)
\nonumber\\
&& \quad
= \left\{ \frac{\partial}{\partial t}
+(\cH_N^{(\beta)}(t)-E_{N,0}^{(\beta)}(t)) \right\}
-\left( \frac{\partial}{\partial t}
-\frac{1}{2} \sum_{j=1}^N \frac{\partial^2}{\partial x_j^2} \right)
\nonumber\\
&& \quad
= V_N^{(\beta)}(t, \x)-E_{N,0}^{(\beta)}(t).
\nonumber
\end{eqnarray}
Therefore, we see
$$
\log \widetilde{M}_N^{\u}(t)
= - \frac{\beta}{2} \Big\{
W_N(t, \b^{\u}(t))-W_N(0, \u) \Big\}
-\int_0^t 
V_N^{(\beta)}(s, \b^{\u}(s)) ds 
+ \int_0^t E_{N, 0}^{(\beta)}(s) ds.
$$
On the other hand, by (\ref{eqn:WN_1}) with (\ref{eqn:cW1}), 
\begin{eqnarray}
&& 
- \frac{\beta}{2} \Big\{
W_N(t, \b^{\u}(t))-W_N(0, \u) \Big\}
\nonumber\\
&& \quad
=\log \left[
\frac{\vartheta_1(\overline{b^u}(t)/2 \pi r; N \tau(t))}{\vartheta_1(\overline{u}/2 \pi r; N \tau(0)}
\prod_{1 \leq j < k \leq N}
\frac{\vartheta_1((b_k^{u_k}(t)-b_j^{u_j}(t))/2 \pi r; N \tau(t))}
{\vartheta_1((u_k-u_j)/2 \pi  r; N \tau(0))}
\right]^{\beta/2}.
\nonumber
\end{eqnarray}
The integral formula (\ref{eqn:Dedekind3}) and the definition (\ref{eqn:hN_1})
prove that (\ref{eqn:MN}) is equal to (\ref{eqn:tildeMN1}). 
Hence $M_N^{\u}(t)$ is martingale, since It\^o's formula gives
$$
dM_N^{\u}(t)=d \widetilde{M}_N^{\u}(t)
= - \frac{\beta}{2} \sum_{j=1}^N \frac{\partial W_N}{\partial x_j}(t, \b^{\u}(t))
M_N^{\u}(t) d b_j^{u_j}(t).
$$
(ii) 
Consider the process
$\widetilde{\b}(t)=(\widetilde{b}_1(t), \dots, \widetilde{b}_N(t))$ solving the 
system of SDEs,
\begin{equation}
d \widetilde{b}_j(t) = db_j^{u_j}(t)
+ \frac{\partial W_N}{\partial x_j}(t, \b^{\u}(t)) dt,
\quad 1 \leq j \leq N, t \in [0, t_*).
\label{eqn:b_tilde1}
\end{equation}
Since
$$
d(\widetilde{b}_j(t) M_N^{\u}(t))
=\widetilde{b}_j(t) d M_N^{\u}(t) + M_N^{\u}(t) d b_j(t),
\quad 1 \leq j \leq N, \quad t \in [0, t_*),
$$
$\widetilde{b}_j(t) M_N^{\u}(t), 1 \leq j \leq N$, are
martingales with respect to $\rP$.
By the factor $(h_N^r(t, \b(t)))^{\beta/2}, \beta >0$ in (\ref{eqn:MN}),
the measure $\P^{\u}$, whose Radon-Nikodym derivative
with respect to $\rP$ is given by (\ref{eqn:RadonN}), 
is concentrated on paths with $\sigma^{\u} > t$.
With respect to $\P^{\u}$,
$\widetilde{\b}(t), t \geq 0$ can be
regarded as the $N$-dimensional BM.
Since (\ref{eqn:b_tilde1}) gives a set of SDEs,
$$
b_j^{u_j}(t)=u_j+\widetilde{b}_j(t)
-\frac{\beta}{2} \int_0^t 
\frac{\partial W_N}{\partial x_j}(s, \b^{\u}(s)) ds,
\quad 1 \leq j \leq N, t \in [0, t_*),
$$
which is equivalent with (\ref{eqn:eDYS1}),
$\b(t), t \in [0, t_*)$ is eDYS$^{(\beta)}$ with respect to
$\P^{\u}$. The proof is hence completed. \qed
\vskip 0.3cm

\SSC{Determinantal Property of eDYS$^{(2)}$
\label{sec:eDYS2}}
\subsection{Observables of eDYS$^{(2)}$ and processes $\Y_{\sp(N)}$} \label{sec:Ykappa}
From now on, we consider eDYS$^{(\beta)}$ with
$$
\beta=2,
$$
in which Lemma \ref{thm:VN} gives $V_N^{(2)}(t, \x) \equiv 0$.
Let $\mM$ be the space of nonnegative integer-valued Radon measures 
on the interval $[0, 2 \pi r)$, 
which is a Polish space with the vague topology.
Any element $\xi$ of $\mM$ can be represented as
$\xi(\cdot) = \sum_{j \in \I}\delta_{x_j}(\cdot)$ 
with an index set $\I$, in which
the sequence of points in $[0, 2 \pi r)$, $\x =(x_j)_{j \in \I}$, 
satisfies $\xi(K)=\sharp\{x_j : x_j \in K \} < \infty$ 
for any subset $K \subset [0, 2 \pi r)$.
Now we consider eDYS$^{(2)}$, $\X(t), t \in [0, t_*)$
as an $\mM$-valued process
and write it as
\begin{equation}
\Xi(t, \cdot)=\sum_{j=1}^N \delta_{X_j(t)}(\cdot),
\quad t \in [0, t_{*}),
\label{eqn:Xi1}
\end{equation}
where $\delta_x(\cdot)$ denotes the Dirac measure concentrated on 
$x \in [0, 2 \pi r)$;
$\delta_x(A)=\1(x \in A)$, $A \in \mB([0, 2 \pi r))$. 
The probability law of $\Xi(t, \cdot), t \in [0, t_{*})$
starting from a fixed configuration $\xi \in \mM$
is denoted by $\P^{\xi}$ 
and the process
specified by the initial configuration
is expressed by
$(\Xi(t), t \in [0,t_{*}), \P^{\xi})$.
The expectations with respect to $\P^{\xi}$
is denoted by $\E^{\xi}$.
We introduce a filtration $\{\cF_{\Xi}(t) : t \in [0,t_{*}) \}$
generated by $\Xi(t), t \in [0, t_*)$, 
which satisfies the usual conditions. 
We set 
\begin{equation}
\mM_0= \{ \xi \in \mM : 
\xi(\{x\}) \leq 1 \mbox { for any }  x \in [0, 2 \pi r) \},
\label{eqn:mM0}
\end{equation}
which denotes a collection of configurations
without any multiple points.

Assume that $\xi=\sum_{j=1}^N \delta_{u_j} \in \mM_0$
with $\u \in \check{\A}_N^{2 \pi r}$.
Let $T \in [0, t_*)$.
Proposition \ref{thm:GirsanovN} implies that
for any $\cF_{\Xi}(T)$-measurable bounded function $F$,
\begin{equation}
\E^{\xi}[F(\Xi(\cdot))]
= \rE^{\u} \left[ 
F\left( \sum_{j=1}^N \delta_{b_j(\cdot)} \right)
\1(\sigma > T) 
\frac{h_N^r(T, \b(T))}{h_N^r(0, \u)} \right],
\label{eqn:ExpN1}
\end{equation}
where 
$\sigma = \inf\{t>0 : \b(t) \in \partial \check{\A}^{2 \pi r}_N\}$.
We assume that 
each $\cF_{\Xi}(T)$-measurable function $F$ 
is given by
\begin{equation}
F(\Xi(\cdot))=\prod_{m=1}^M g_m(\X(t_m))
\label{eqn:F_gN}
\end{equation}
for an arbitrary $M \in \N$ and for
an arbitrary increasing series of times,
$0 \leq t_1 < \dots < t_M \leq T < t_*$ 
with bounded measurable functions $g_m, 1 \leq m \leq M$.

For $N$ variables $\x=(x_j)_{j=1}^N$, elementary symmetric polynomials
are defined as
\begin{equation}
e_r(\x)=\sum_{1 \leq j_1 < j_2 < \cdots < j_r \leq N}
x_{j_1} x_{j_2} \cdots x_{j_r},
\quad r=1,2, \dots, N.
\label{eqn:sym1}
\end{equation}
If $g(\x)$ is a symmetric function of $\x$,
it can be regarded as a function of $(e_r(\x))_{r=1}^N$,
which we write as
$g(\x)=\widetilde{g}(e_1(\x), e_2(\x), \dots, e_N(\x))$.

We define the observables for eDYS$^{(2)}$ considered below.

\begin{df}
\label{thm:observableN}
For $T \in [0, t_*)$, 
define the $\cF_{\Xi}(T)$-observables of eDYS$^{(2)}$ 
as a set of $\cF_{\Xi}(T)$-measurable bounded functions 
given in the form (\ref{eqn:F_gN}) satisfying the following conditions. \\
{\rm (1)} 
$\{g_m(\x) \}_{m=1}^M$ are periodic functions, 
\begin{equation}
g_m((x_j+2 \pi r n_j)_{j=1}^N) =g_m( \x), \quad \forall n_j \in \Z, \quad 1 \leq j \leq N, \quad 1 \leq m \leq M.
\label{eqn:g_periodN}
\end{equation}
{\rm (2)} 
$\{g_m(\x) \}_{m=1}^M$ are symmetric functions of variables $\x=(x_j)_{j=1}^N$.
That is, 
\begin{equation}
g_m(\pi(\x))=g_m(\x), \quad
\forall \pi \in \mS_N, \quad 1 \leq m \leq M.
\label{eqn:g_symN}
\end{equation}
{\rm (3)}
$\{g_m(\x)\}_{m=1}^M$ are invariant under reflections of $\overline{x}$ 
at the boundary $\{0, 2 \pi r\}$ of the interval $(0, 2 \pi r)$.
Since $e_1(\x)=\kappa_N+\overline{x}$, this condition is written as
\begin{eqnarray}
\widetilde{g}(\kappa_N-x, (e_r(\x))_{r=2}^N)
&=& \widetilde{g}(\kappa_N+x, (e_r(\x))_{r=2}^N),
\nonumber\\
\widetilde{g}(\kappa_N+ 2 \pi r +x, (e_r(\x))_{r=2}^N)
&=& \widetilde{g}(\kappa_N+ 2 \pi r -x, (e_r(\x))_{r=2}^N),
\quad x \in (0, 2 \pi r).
\label{eqn:e_1}
\end{eqnarray}
\end{df}

Note that $h_N^r(t, \x)$ is anti-symmetric 
for any exchange of two particle positions
as shown by (\ref{eqn:hN_3}).
The expression (\ref{eqn:hN_1}) implies that
$h_N^r(t, \x)$ changes its sign for reflection
of $\overline{x}$ at the boundary $\{0, 2 \pi r\}$ of the
interval $(0, 2 \pi r)$, since
$\vartheta_1(-\overline{x}/2 \pi r; N \tau(t))=
-\vartheta_1(\overline{x}/2 \pi r; N \tau(t))$
and
$\vartheta_1((2 \pi r+\overline{x})/2 \pi r; N \tau(t))
=-\vartheta_1((2 \pi r-\overline{x})/2 \pi r; N \tau(t))$.
Therefore, by the reflection principle of BM
and the fact that
$\overline{\check{X}}(t)$ of eDYS$^{(2)}$
is the time change of eBES$^{(3)}$, 
$\check{X}(t)$, 
(see Sec. \ref{sec:eDYS_beta}),
if $F$ is an $\cF_{\Xi}(T)$-observable, 
we can omit the condition $\1(\sigma > T)$ in (\ref{eqn:ExpN1}),
\begin{equation}
\E^{\xi}[F(\Xi(\cdot))]
= \rE^{\u} \left[ 
F\left( \sum_{j=1}^N \delta_{b_j(\cdot)} \right)
\frac{h_N^r(T, \b(T))}{h_N^r(0, \u)} \right],
\label{eqn:ExpN2}
\end{equation}
since all paths with $\sigma < T$ are canceled out.

If we use the determinantal expression (\ref{eqn:hN_2})
of $h^r_N$ with $\beta=2$, (\ref{eqn:ExpN2}) is written as
\begin{eqnarray}
\label{eqn:ExpN3}
\E^{\xi}[F(\Xi(\cdot))]
&=& \rE^{\u} \left[
F\left( \sum_{j=1}^N \delta_{b_j(\cdot)} \right)
\frac{q^r_N(t_*-T, \v|\b(T))}{q^r_N(t_*, \v|\u)} \right]
\\
\label{eqn:ExpN4}
&=& \sum_{o \in\Z} \mu_N^{\u}(o)
\rE^{\u} \left[
F\left( \sum_{j=1}^N \delta_{b_j(\cdot)} \right)
\frac{q_{\W_N}(t_*-T, \v_{o}|\b(T))}{q_{\W_N}(t_*, \v_o|\u)} \right],
\end{eqnarray}
with
\begin{equation}
\mu_N^{\u}(o)=\frac{(-1)^o q_{\W_N}(t_*, \v_o|\u)}
{\sum_{o' \in \Z} (-1)^{o'} q_{\W_N}(t_*, \v_{o'}|\u)}
\label{eqn:mu_N1}
\end{equation}
where $q_{\W_N}$ is the Karlin-McGregor determinant 
of $p^{\rm BM}$ given by (\ref{eqn:KM1}), and
$$
\v_o=(v_{j+o})_{j=1}^N \quad
\mbox{with} \quad
v_{j+o}=\frac{2 \pi r}{N}(j+o-1), \quad 1 \leq j \leq N.
$$
The integer $o$ is called an offset in \cite{LW13},
which is given by $o=N w+\ell$ for the $N$ particle systems
with the winding number $w$ (see Eq.(\ref{eqn:pr+1}))
and with the shift by $\ell \in \{0,1, \dots, N-1\}$
(see Eq.(\ref{eqn:detp+1})).
The expression (\ref{eqn:ExpN4}) shows that
for each $o \in \Z$ we consider an $N$-dimensional
Brownian bridge of time duration $t_*$, 
$\b(t)=(b_1(t), \dots, b_N(t)), t \in [0, t_*)$, which is started at 
$\u$ and ended at $\v_{o}$,
and take the summation over $o \in \Z$ with the
signed measure (\ref{eqn:mu_N1}).

Assume that $F$ is given by (\ref{eqn:F_gN}).
Then (\ref{eqn:ExpN3}) is written as
\begin{eqnarray}
\E^{\xi}[F(\Xi(\cdot))]
&=&\prod_{m=1}^M \int_{\R^N} d \x^{(m)}
\prod_{m=1}^M \Big\{ g_m(\x^{(m)}) p^{\rm BM}(t_m-t_{m-1}, \x^{(m)}|\x^{(m-1)}) \Big\}
\nonumber\\
&& \qquad \qquad \qquad \times
\frac{q^r_N(t_*-t_M, \v|\x^{(M)})}{q^r_N(t_*, \v|\u)},
\label{eqn:ExpN5}
\end{eqnarray}
where $d \x^{(m)}=\prod_{j=1}^N dx^{(m)}_j$ and
\begin{equation}
p^{\rm BM}_N(t, \y|\x)=
\prod_{j=1}^N p^{\rm BM}(t, y_j|x_j),
\quad \x, \y \in \R^N.
\label{eqn:pBM_N}
\end{equation}

\begin{df}
\label{thm:process_Ykappa}
For $N \in \N$, the process 
$\Y_{\sp(N)}(t)=(Y_{\sp(N), 1}(t), \dots, Y_{\sp(N), N}(t))$, $t \in [0, t_*)$,
is defined in $\A^{2 \pi r}_N$
by the following finite-dimensional distributions.
Assume that the initial state is given by $\u \in \check{\A}^{2 \pi r}_N$.
For an arbitrary $M \in \N$, 
an arbitrary strictly increasing sequence of times
$0 \leq t_1 < \cdots < t_M < t_*$,
and arbitrary Borel sets $A_m \in \mB(\A^{2 \pi r}_N)$, $1 \leq m \leq M$, 
\begin{eqnarray}
&& \bP^{\u} [\Y_{\sp(N)}(t_1) \in A_1, \Y_{\sp(N)}(t_2) \in A_2, \dots,
\Y_{\sp(N)}(t_M) \in A_M]
\nonumber\\
&& \quad =
\prod_{m=1}^M \int_{A_m} d \x^{(m)}
\prod_{m=1}^M q^r_N(t_m-t_{m-1}, \x^{(m)}|\x^{(m-1)}),
\label{eqn:bE_N1}
\end{eqnarray}
where $t_0=0, \x^{(0)}=\u$.
The expectation is denoted by $\bE^{\u}$.
The natural filtration of $\Y_{\sp(N)}$ is denoted by
$\cF_{Y_{\sp(N)}}(t)=\sigma(\Y_{\sp(N)}(s), s \in [0, t])$
for $t \in [0, t_*)$.
\end{df}

We will prove the following.
\begin{lem}
\label{thm:q-formN1}
Let $T \in [0, t_*)$. For any $\cF_{\Xi}(T)$-observable $F$,
\begin{eqnarray}
\E^{\xi}[F(\Xi(\cdot))]
&=& \bE^{\u} \left[ F\left( \sum_{j=1}^N \delta_{Y_{\sp(N), j}(\cdot)} \right)
\frac{q^r_N(t_*-T, \v | \Y_{\sp(N)}(T))}{q^r_N(t_*, \v| \u)} \right], 
\label{eqn:q-formN1}
\end{eqnarray}
$\xi = \sum_{j=1}^N \delta_{u_j}$, 
$\u \in \check{\A}^{2 \pi r}_N$.
\end{lem}
\noindent{\it Proof.} \,
By the definition (\ref{eqn:pkappa}), for $\ell \in \Z$,
\begin{eqnarray}
p^r_{\sp(N)}(t, y|x+2 \pi r \ell)
&=& \sum_{w \in \Z} \sp(N)^w p^{\rm BM}(t, y+2 \pi r w|x+2 \pi r \ell)
\nonumber\\
&=& \sp(N)^{\ell} \sum_{w' \in \Z}
\sp(N)^{w'} p^{\rm BM}(t, y+2 \pi r w'|x)
\nonumber\\
&=& \sp(N)^{\ell} p^r_{\sp(N)}(t, y|x),
\label{eqn:pN_1}
\end{eqnarray}
and hence, for $\ell_j \in \Z, 1 \leq j \leq N$,
\begin{equation}
q^r_N(t, \y|(x_j+2 \pi r \ell_j)_{j=1}^N)
=\prod_{j=1}^N \sp(N)^{\ell_j} q^r_N(t, \y|\x).
\label{eqn:pN_2}
\end{equation}
We see 
\begin{eqnarray}
&& \int_{\R^N} d \x^{(M)} g_M(\x^{(M)})
p^{\rm BM}_N(t_M-t_{M-1}, \x^{(M)}|\x^{(M-1)})
q^r_N(t_*-t_M, \v|\x^{(M)})
\nonumber\\
&& =
\prod_{j=1}^N \sum_{\ell_j \in \Z}
\int_{2 \pi r \ell_j}^{2 \pi r (\ell_j+1)} dx^{(M)}_j \,
g_M(\x^{(M)}) p^{\rm BM}_N(t_M-t_{M-1}, \x^{(M)}|\x^{(M-1)})
q^r_N(t_*-t_M, \v|\x^{(M)})
\nonumber\\
&& =
\prod_{j=1}^N \sum_{\ell_j \in \Z}
\int_{0}^{2 \pi r} dy^{(M)}_j \, g_M((y^{(M)}_j+2 \pi r \ell_j)_{j=1}^N)
p^{\rm BM}_N(t_M-t_{M-1}, (y^{(M)}_k+2 \pi r \ell_k)_{k=1}^N | \x^{(M-1)})
\nonumber\\
&& \qquad \qquad \qquad \qquad \times
q^r_N(t_*-t_M, \v | (y^{(M)}_n+2 \pi r \ell_n)_{n=1}^N),
\nonumber
\end{eqnarray}
where we have put $y^{(M)}_j=x^{(M)}_j-2 \pi r \ell_j$
for each $1 \leq j \leq N$.
Since $g_M$ is assumed to be periodic in the sense (\ref{eqn:g_periodN}),
the above is written using (\ref{eqn:pN_2}) as
\begin{eqnarray}
&& \int_{[0, 2 \pi r)^N} d \y^{(M)} \,
g_M(\y^{(M)}) 
\prod_{j=1}^N \left\{
\sum_{\ell_j \in \Z} \sp(N)^{\ell_j}
p^{\rm BM}(t_M-t_{M-1}, y^{(M)}_j+2 \pi r \ell_j | x^{(M-1)}_j) \right\}
\nonumber\\
&& \qquad \qquad \qquad \qquad \times
q^r_N(t_*-t_M, \v|\y^{(M)})
\nonumber\\
&& = \int_{[0, 2 \pi r)^N} d \y^{(M)} \,
g_M(\y^{(M)}) 
\prod_{j=1}^N 
p^r_{\sp(N)}(t_M-t_{M-1}, y^{(M)}_j | x^{(M-1)}_j)
q^r_N(t_*-t_M, \v|\y^{(M)}). 
\nonumber
\end{eqnarray}
If we write
$\pi(\A^{2 \pi r}_N)
=\{\x \in \R^N:  x_{\pi(1)} < \cdots
< x_{\pi(N)} < x_{\pi(1)}+2 \pi r \}$
for $\pi \in \mS_N$, the above is further rewritten as
\begin{eqnarray}
&& \sum_{\pi \in \mS_N} \int_{\pi(\A^{2 \pi r}_N)} d \y^{(M)} \,
g_M(\y^{(M)})
\prod_{j=1}^N p^r_{\sp(N)}(t_M-t_{M-1}, y^{(M)}_j|x^{(M-1)}_j)
q^r_N(t_*-t_M, \v|\y^{(M)})
\nonumber\\
&& = \sum_{\pi \in \mS_N} \int_{\A^{2 \pi r}_N} d \y^{(M)} \,
g_M(\pi(\y^{(M)}))
\prod_{j=1}^N p^r_{\sp(N)}(t_M-t_{M-1}, y^{(M)}_{\pi(j)}|x^{(M-1)}_j)
q^r_N(t_*-t_M, \v| \pi(\y^{(M)}))
\nonumber\\
&& = \int_{\A^{2 \pi r}_N} d \y^{(M)} \,
g_M(\y^{(M)}) \sum_{\pi \in \mS_N}  {\rm sgn}(\pi)
\prod_{j=1}^N p^r_{\sp(N)}(t_M-t_{M-1}, y^{(M)}_{\pi(j)}|x^{(M-1)}_j)
\nonumber\\
&& \qquad \qquad \qquad \qquad \times
q^r_N(t_*-t_M, \v| \y^{(M)}), 
\nonumber
\end{eqnarray}
where the assumption that $g_M$ is symmetric in the sense of (\ref{eqn:g_symN})
and the basic property of determinant, 
$q^r_N(\cdot, \v|\pi(\y))={\rm sgn}(\pi) q^r_N(\cdot, \v|\y)$, $\pi \in \mS_N$, 
were used.
Then we have the equality
\begin{eqnarray}
&& \int_{\R^N} d \x^{(M)} \, g_M(\x^{(M)})
p^{\rm BM}_N(t_M-t_{M-1}, \x^{(M)}| \x^{(M-1)}) q^r_N(t_*-t_M, \v|\x^{(M)})
\nonumber\\
&& = \int_{\A^{2 \pi r}_N} d \x^{(M)} \, g_M(\x^{(M)})
q^r_N(t_M-t_{M-1}, \x^{(M)}| \x^{(M-1)}) q^r_N(t_*-t_M, \v|\x^{(M)}).
\nonumber
\end{eqnarray}
Similarly, we can prove the following equalities for $1 \leq m \leq M-1$,
\begin{eqnarray}
&& \int_{\R^N} d \x^{(m)} \, g_m(\x^{(m)})
p^{\rm BM}_N(t_m-t_{m-1}, \x^{(m)}| \x^{(m-1)}) q^r_N(t_{m+1}-t_m, \x^{(m+1)} | \x^{(m)})
\nonumber\\
&& = \int_{\A^{2 \pi r}_N} d \x^{(m)} \, g_m(\x^{(m)})
q^r_N(t_m-t_{m-1}, \x^{(m)}| \x^{(m-1)}) q^r_N(t_{m+1}-t_m, \x^{(m+1)} |\x^{(m)}).
\nonumber
\end{eqnarray}
Combining these equalities, the equivalence between (\ref{eqn:ExpN5}) 
and (\ref{eqn:q-formN1}) is proved. \qed

\subsection{Spatio-temporal correlations of eDYS$^{(2)}$} \label{sec:determinantal}

Comparing the product expression (\ref{eqn:hN_1}) with the
determinantal expression (\ref{eqn:hN_2}) of $h^r_N(t, \x)$, we see
\begin{eqnarray}
&& \frac{q^r_N(t_*-t, \v|\x)}{q^r_N(t_*, \v|\u)}
=\left( \frac{\eta(N\tau(t))}{\eta(N\tau(0))} \right)^{-(N-1)(N-2)/2}
\nonumber\\
&& \qquad \qquad \times
\frac{\vartheta_1(\overline{x}/2 \pi r; N \tau(t))}{\vartheta_1(\overline{u}/2 \pi r; N \tau(0))}
\prod_{1 \leq j < k \leq N}
\frac{\vartheta_1((x_k-x_j)/2 \pi r; N \tau(t))}{\vartheta_1((u_k-u_j)/2 \pi r; N \tau(0))},
\label{eqn:MDR1}
\end{eqnarray}
$t \in [0,t_*), \u, \x \in \check{\A}^{2 \pi r}_N$.

Let $\widetilde{b}_j(t), t \geq 0, 1 \leq j \leq N$ be
independent one-dimensional standard BMs started at the origin.
We write the expectation with respect to 
$\widetilde{\b}(t)=(\widetilde{b}_1(t), \dots, \widetilde{b}_N(t)), t \geq 0$
as $\widetilde{\rE}[\, \cdot \,]$.
In the previous paper \cite{Kat15},
we showed that, using the determinantal equalities
involving Jacobi's theta functions proved in
the paper by Rosengren and Schlosser \cite{RS06}
(see also \cite{TV97,War02,KN03,Kra05}),
(\ref{eqn:MDR1}) is expressed by the following
(Lemmas 2.11 and 2.12 in \cite{Kat15}),
\begin{eqnarray}
\frac{q^r_N(t_*-t, \v|\x)}{q^r_N(t_*, \v|\u)}
&=& \widetilde{\rE} \left[
\frac{q^r_N(t_*, \v| (x_j+i \widetilde{b}_j(t))_{j=1}^N)}{q^r_N(t_*, \v|\u)}
\right]
\nonumber\\
&=& \widetilde{\rE} \left[
\det_{1 \leq j, k \leq N}
\left[ \Phi_{\xi}^{u_k}(x_j+ i \widetilde{b}_j(t)) \right] \right]
\nonumber\\
&=& \int_{\R^N} d \widetilde{\b} \,
p^{\rm BM}_N(t, \widetilde{\b}| \0)
\det_{1 \leq j, k \leq N}
\left[ \Phi_{\xi}^{u_k}(x_j+ i \widetilde{b}_j) \right],
\label{eqn:MDR2}
\end{eqnarray}
$t \in [0, t_*), \u, \x \in \check{\A}^{2 \pi r}_N$
with
\begin{equation}
\Phi_{\xi}^{u_k}(z)
= \frac{\vartheta_1((\overline{u}+z-u_k)/2 \pi r; N \tau(0))}
{\vartheta_1(\overline{u}/2 \pi r; N \tau(0))}
\prod_{\substack{1 \leq \ell \leq N, \cr \ell \not=k}}
\frac{\vartheta_1((z-u_{\ell})/2 \pi r; N \tau(0))}
{\vartheta_1((u_k-u_{\ell})/2 \pi r; N \tau(0))}, \, 1 \leq k \leq N.
\label{eqn:Phi1}
\end{equation}
Since (\ref{eqn:setting}) gives $\tau(0)=i t_*/2 \pi r^2$,
and thus, $\Im \tau(0) > 0$,
$\Phi_{\xi}^{u_k}(z), 1 \leq k \leq N$ are holomorphic for $|z| < \infty$.

By Jacobi's imaginary transformation (\ref{eqn:Jacobi_imaginary}), 
we obtain the equality
$$
\vartheta_1\left( \frac{x}{2 \pi r}; N \tau(0) \right)
= i \sqrt{\frac{2 \pi r^2}{Nt_*}} e^{-x^2/2Nt_*}
\vartheta_1 \left( - \frac{i r x}{N t_*}; \frac{2 \pi i r^2}{N t_*} \right),
$$
where $\tau(0)=i t_*/2 \pi r^2$ was used, and 
we can rewrite (\ref{eqn:Phi1}) as 
\begin{eqnarray}
\Phi_{\xi}^{u_k}(z)
&=& \exp \left[
-\frac{1}{2 t_*} \left\{ z \left(z - \frac{2 \kappa_N}{N} \right)
-u_k \left( u_k -\frac{2 \kappa_N}{N} \right)\right\} \right]
\nonumber\\
&\times& \frac{\vartheta_1(-ir(\overline{u}+z-u_k)/Nt_*; 2 \pi i r^2/N t_*)}
{\vartheta_1(-ir \overline{u}/N t_*; 2 \pi i r^2/N t_*)}
\nonumber\\
&\times& 
\prod_{\substack{1 \leq \ell \leq N, \cr \ell \not=k}}
\frac{\vartheta_1(-ir(z-u_{\ell})/Nt_*; 2 \pi i r^2/N t_*)}
{\vartheta_1(-ir (u_k-u_{\ell}) /N t_*; 2 \pi i r^2/N t_*)}, 
\quad 1 \leq k \leq N.
\label{eqn:Phi2}
\end{eqnarray}
For $t \in [0, t_*)$ define
\begin{eqnarray}
\cM_{\xi}^{u_k}(t, x)
&=& \widetilde{\rE}
\Big[ \Phi_{\xi}^{u_k}(x+i \widetilde{b}(t)) \Big]
\nonumber\\
&=& \int_{\R} d \widetilde{b} \
\frac{e^{-\widetilde{b}^2/2t}}{\sqrt{2 \pi t}}
\Phi_{\xi}^{u_k}(x+i \widetilde{b}),
\quad 1\leq k \leq N.
\label{eqn:cM1}
\end{eqnarray}
We prove the following.
\begin{lem}
\label{thm:martingale}
Assume $\u \in \check{\A}^{2 \pi r}_N$ 
and $\xi=\sum_{j=1}^N \delta_{u_j} \in \mM_0$.
Then the following are satisfied. 
\begin{description}
\item{\rm (i)} \,
$\cM_{\xi}^{u_k}(t, Y_{\sp(N)}(t))$, $1 \leq k \leq N$, $t \in [0, t_*)$
are martingales in the sense that
\begin{equation}
\bE[ \cM_{\xi}^{u_k}(t, Y_{\sp(N)}(t)) | \cF_{Y_{\sp(N)}}(s)]
=\cM_{\xi}^{u_k}(s, Y_{\sp(N)}(s))
\quad \mbox{a.s.}
\label{eqn:martingale_1}
\end{equation}
for any $0 \leq s \leq t < t_*$.

\item{\rm (ii)} \,
For any $t \in [0, t_*)$, 
$\cM_{\xi}^{u_k}(t, x)$, $1 \leq k \leq N$,
are linearly independent functions of $x \in [0, 2 \pi r)$.

\item{\rm (iii)} \,
For $1 \leq j, k \leq N$,
$$
\cM_{\xi}^{u_k}(0,u_j)=\delta_{jk}.
$$
\end{description}
\end{lem}
\noindent{\it Proof.} \,
(i) For the quasi-periodicity (\ref{eqn:quasi_periodic}) of $\vartheta_1$,
the definition (\ref{eqn:cM1}) with the expression (\ref{eqn:Phi1})
implies that, for $w \in \Z$,
\begin{equation}
\cM_{\xi}^{u_k}(t, x+2 \pi r w)
= \sp(N)^w \cM_{\xi}^{u_k}(t, x),
\quad 1 \leq k \leq N.
\label{eqn:M_quasi}
\end{equation}
On the other hand, (\ref{eqn:pkappa}) gives
\begin{eqnarray}
&& \bE[ \cM_{\xi}^{u_k}(t, Y_{\sp(N)}(t)) | \cF_{Y_{\sp(N)}}(s)]
\nonumber\\
&& \quad = \int_0^{2 \pi r} dy \,
\cM_{\xi}^{u_k}(t, y) p^r_{\sp(N)}(t-s, y|Y_{\sp(N)}(s))
\nonumber\\
&& \quad = \sum_{w \in \Z} \int_{2 \pi r w}^{2 \pi r (w+1)} dy \,
\cM_{\xi}^{u_k}(t, y-2 \pi r w) \sp(N)^w p^{\rm BM}(t-s, y|Y_{\sp(N)}(s)),
\quad 1\leq k \leq N.
\nonumber
\end{eqnarray}
By (\ref{eqn:M_quasi}), this is equal to
$$
\int_{\R} dy \, \cM_{\xi}^{u_k}(t, y)
p^{\rm BM}(t-s, y|Y_{\sp(N)}(s)) \quad \mbox{a.s.}
$$
Next we consider the expression (\ref{eqn:Phi2}).
By the definition of $\vartheta_1$ given in (\ref{eqn:theta}), we
obtain the following expansions,
\begin{eqnarray}
&&
\frac{\vartheta_1(-ir(\overline{u}+z-u_k)/Nt_*; 2 \pi i r^2/N t_*)}
{\vartheta_1(-ir \overline{u}/N t_*; 2 \pi i r^2/N t_*)}
=\sum_{n_0 \in \Z} c^0_{n_0} e^{(2n_0-1) \pi r z/N t_*},
\nonumber\\
&&
\frac{\vartheta_1(-ir(z-u_{\ell})/Nt_*; 2 \pi i r^2/N t_*)}
{\vartheta_1(-ir (u_k-u_{\ell}) /N t_*; 2 \pi i r^2/N t_*)}
= \sum_{n_{\ell} \in \Z} c^{\ell}_{n_{\ell}} e^{(2n_{\ell}-1) \pi r z/N t_*},
\quad 1 \leq \ell \leq N, \ell \not= k,
\nonumber
\end{eqnarray}
with
\begin{eqnarray}
c^0_{n_0} &=&
\frac{i (-1)^{n_0}}{\vartheta_1(-ir \overline{u}/N t_*; 2 \pi i r^2/N t_*)}
\exp \left[ - \frac{2 \pi^2 r^2}{Nt_*} \left( n_0-\frac{1}{2} \right)^2
-(2n_0-1) \frac{ir}{N t_*}(\overline{u}-u_k) \right]
\nonumber\\
c^{\ell}_{n_{\ell}} &=&
\frac{i (-1)^{n_{\ell}}}{\vartheta_1(-ir (u_k-u_{\ell}) /N t_*; 2 \pi i r^2/N t_*)}
\exp \left[ - \frac{2 \pi^2 r^2}{Nt_*} \left( n_{\ell}-\frac{1}{2} \right)^2
-(2n_{\ell}-1) \frac{ir}{N t_*}u_{\ell} \right],
\nonumber\\
&& \hskip 5cm 1 \leq \ell \leq N, \quad \ell \not=k.
\nonumber
\end{eqnarray}
For each $1 \leq k \leq N$, 
we introduce an $N$-component index
$\n =(n_0, n_1, \dots, n_{k-1}, n_{k+1}, \dots, n_N) \in \Z^N$,
and put
$$
c_{\n}^{(k)}= \exp
\left[\frac{1}{2 t_*} u_k \left( u_k-\frac{2 \kappa_N}{N} \right) \right]
\prod_{\substack{ 0 \leq \ell \leq N, \cr \ell \not= k}}
c_{n_{\ell}}^{\ell}, \quad
1 \leq k \leq N,
$$
with (\ref{eqn:kappaN}). 
Then we obtain the following expression for the
entire functions (\ref{eqn:Phi2}),
$$
\Phi_{\xi}^{u_k}(z)
= \exp \left[ - \frac{1}{2t_*} z \left( z-\frac{2 \kappa_N}{N} \right) \right]
\sum_{\n \in \Z^N} c_{\n}^{(k)} \exp \left[
\frac{\pi r z}{N t_*} \sum_{\substack{ 0 \leq \ell \leq N, \cr \ell \not=k}}
(2 n_{k}-1) \right],
\quad 1 \leq k \leq N.
$$
By the definition (\ref{eqn:cM1}), we have the following expansion
formula,
\begin{eqnarray}
&& \cM_{\xi}^{u_k}(t, x)
= \exp \left[ - \frac{1}{2t_*} x \left( x-\frac{2 \kappa_N}{N} \right) \right]
\sum_{\n \in \Z^N} c_{\n}^{(k)} \exp \left[
\frac{\pi r x}{N t_*} \sum_{\substack{ 0 \leq \ell \leq N, \cr \ell \not=k}}
(2 n_{k}-1) \right]
\nonumber\\
&& \quad \times
\widetilde{\rE} \left[
\exp \left[ \frac{\widetilde{b}(t)^2}{2t_*} 
-i \frac{\widetilde{b}(t)}{t_*} 
\left( x - \frac{\kappa_N}{N} 
-\frac{\pi r}{2N}
\sum_{\substack{0 \leq \ell \leq N, \cr \ell \not= k}}
(2 n_{\ell}-1) \right)
\right]
\right], \quad 1 \leq k \leq N.
\nonumber
\end{eqnarray}
They are calculated as
\begin{equation}
\cM_{\xi}^{u_k}(t, x)
= \sqrt{2 \pi t_*} e^{\alpha^2/2t_*} 
\sum_{\n \in \Z^N} c_{\n}^{(k)}
p^{\rm BM} \left( t_*-t, \left. 
\frac{\kappa_N}{N}+\frac{\pi r}{N} 
\sum_{\substack{0 \leq \ell \leq N, \cr \ell \not=k}}
(2 n_{\ell}-1) \right| x \right),
\quad 1 \leq k \leq N,
\label{eqn:M_series2}
\end{equation}
where $p^{\rm BM}$ denotes the tpd of a single BM.
The Chapman-Kolmogorov equation 
$$
\int_{\R} dy \, p^{\rm BM}(t_*-t, \alpha|y)
p^{\rm BM}(t-s, y|x)=p^{\rm BM}(t_*-s, \alpha|x),
\quad 0 < s < t < t_*, \quad \alpha \in \R, 
$$
implies that $p^{\rm BM}(t_*-t, \alpha| b(t))$, $\alpha \in \R$ is
a continuous martingale with respect to the natural filtration 
generated by BM, 
$\cF_{b}(t)=\sigma(b(s): s \in [0, t])$ in the time period $t \in [0, t_*)$;
$$
\rE[p^{\rm BM}(t_*-t, \alpha|b(t)) | \cF_{b}(s)]
=p^{\rm BM}(t_*-s, \alpha|b(s)) \quad \mbox{a.s.}
$$
for any two stopping times $0 \leq s \leq t < t_*$.
Since $p(t_*-t, \alpha|x)$ gives the Gaussian distribution
of $\alpha$ with mean $x$ and variance $t_*-t \in (0, t_*]$ for $t \in [0, t_*)$, 
convergence of the series (\ref{eqn:M_series2})
is obvious.
Thus (\ref{eqn:martingale_1}) is proved.
By assumption $\xi \in \mM_0$, and thus
the zeroes of $\Phi_{\xi}^{u_j}(z)$ are distinct
from those of $\Phi_{\xi}^{u_k}(z)$, if $j \not= k$.
Hence (ii) is proved.
By (\ref{eqn:cM1}), 
$\cM_{\xi}^{u_k}(0, x)=\lim_{t \downarrow 0}
\widetilde{\rE}[\Phi_{\xi}^{u_k}(x+i \widetilde{b}(t))]
=\Phi_{\xi}^{u_k}(x), 1 \leq k \leq N$.
Since $\Phi_{\xi}^{u_k}(u_j)=\delta_{jk}, 1 \leq j, k \leq N$
by the definition (\ref{eqn:Phi1}),
(iii) is also satisfied. \qed
\vskip 0.3cm

For $\x=(x_j)_{j=1}^N$, let
\begin{equation}
\cD^{\xi}(t, \x)
=\det_{1 \leq j, k \leq N}
[ \cM_{\xi}^{u_k}(t, \x) ],
\quad t \in [0, t_*),
\label{eqn:cD1}
\end{equation}
where $\xi  = \sum_{j=1}^N u_j \in \mM_0$. 
For the process $\Y_{\sp(N)}(t), t \geq 0$ started at
$\u \in \check{\A}^{2 \pi r}_N$,
whose finite-dimensional distributions are given by
(\ref{eqn:bE_N1}), 
consider $\cD^{\xi}(t, \Y_{\sp(N)}(t)), t \in [0, t_*)$. 
By part (i) of Lemma \ref{thm:martingale},
$\cD^{\xi}(t, \Y_{\sp(N)}), t \in [0, t_*)$,
is a martingale, 
and it is not identically zero by part (ii) of Lemma \ref{thm:martingale}.
We call this a {\it determinantal martingale}
and the following the 
{\it determinantal martingale representation} (DMR)
\cite{KT13,Kat14,Kat16_Springer}.

\begin{prop}
\label{thm:DMR}
Let $T \in [0, t_*)$. For any $\cF_{\Xi}(T)$-observable $F$,
the process $(\Xi(t), t \in [0, t_*), \P^{\xi})$,
$\xi \in \mM_0$ has the DMR with respect to
the process $\Y_{\sp(N)}, t \geq 0$, that is, 
\begin{eqnarray}
\E^{\xi}[F(\Xi(\cdot))]
&=& \bE^{\u} \left[ F\left( \sum_{j=1}^N \delta_{Y_{\sp(N), j}(\cdot)} \right)
\cD^{\xi}(t, \Y_{\sp(N)}(T)) \right].
\label{eqn:DMR1}
\end{eqnarray}
\end{prop}
\noindent{\it Proof.} \,
By the multi-linearity of determinant and the independence of
$\widetilde{b}_j(t), 1 \leq j \leq N$,
(\ref{eqn:MDR2}) is equal to (\ref{eqn:cD1}).
Thus (\ref{eqn:q-formN1}) is equal to (\ref{eqn:DMR1}). \qed
\vskip 0.3cm

For any integer $M \in \N$,
a sequence of times
$\t=(t_1,\dots,t_M)$ with 
$0 \leq t_1 < \cdots < t_M \leq T < t_{*}$,
and a sequence of $\cF_{\Xi}(T)$-observables 
$\f=(f_{t_1},\dots,f_{t_M})$,
the {\it moment generating function} of multitime distribution
of $(\Xi(t), t \in [0,t_{*}), \P^{\xi})$ is defined by
\begin{equation}
\Psi^{\xi}_{\t}[\f]
= \E^{\xi} \left[ \exp \left\{ \sum_{m=1}^{M} 
\int_{0}^{2 \pi r} f_{t_m}(x) \Xi(t_m, dx) \right\} \right].
\label{eqn:GF1}
\end{equation}
It is expanded with respect to `test functions' 
$\chi_{t_m}(\cdot)=e^{f_{t_m}(\cdot)}-1, 
1 \leq m \leq M$
as
$$
\Psi^{\xi}_{\t}[\f]
=\sum_
{\substack
{0 \leq N_m \leq N, \\ 1 \leq m \leq M} }
\int_{\prod_{m=1}^{M} \A^{2 \pi r}_{N_m}}
\prod_{m=1}^{M} \left\{ d \x_{N_m}^{(m)}
\prod_{j=1}^{N_{m}} 
\chi_{t_m} \Big(x_{j}^{(m)} \Big) \right\}
\rho^{\xi} 
\Big( t_{1}, \x^{(1)}_{N_1}; \dots ; t_{M}, \x^{(M)}_{N_M} \Big),
$$
and it defines the {\it spatio-temporal correlation functions}
$\rho^{\xi}(\cdot)$ for the process $(\Xi(t), t \in [0,t_{*}), \P^{\xi})$.

Given an integral kernel
$
\mbK(s,x;t,y), 
(s,x), (t,y) \in [0, t_{*}) \times [0, 2 \pi r),
$
the {\it Fredholm determinant} is defined as
\begin{eqnarray}
&& \mathop{{\rm Det}}_
{\substack{
(s,t)\in \{t_1, \dots, t_M\}^2, \\
(x,y)\in [0, 2 \pi r)^2}
}
 \Big[\delta_{st} \delta (x-y)
+ \mbK(s,x;t,y) \chi_{t}(y) \Big]
\nonumber\\
&& 
=\sum_
{\substack
{0 \leq N_m \leq N, \\ 1 \leq m \leq M} }
\sum_
{\substack
{\x^{(m)}_{N_m} \in \check{\A}^{2 \pi r}_{N_m}, 
\\ 1 \leq m \leq M} }
\prod_{m=1}^{M}
\prod_{j=1}^{N_{m}} 
\chi_{t_m} \Big(x_{j}^{(m)} \Big)
\det_{\substack
{1 \leq j \leq N_{m}, 1 \leq k \leq N_{n}, \\
1 \leq m, n \leq M}
}
\Bigg[
\mbK(t_m, x_{j}^{(m)}; t_n, x_{k}^{(n)} )
\Bigg].
\nonumber\\
\label{eqn:F_det}
\end{eqnarray}
We put the following definition \cite{BR05,KT10,Kat16_Springer}.
\begin{df}
\label{thm:determinantal}
For a given initial configuration $\xi$, 
if any moment generating function (\ref{eqn:GF1}) for observables
is expressed by a Fredholm determinant with an integral kernel $\mbK^{\xi}$, we say
the process $(\Xi(t), t \in [0, t_{*}), \P^{\xi})$ is determinantal
for the observables. 
In this case, 
all spatio-temporal correlation functions for observables
are given by determinants as 
$$
\rho^{\xi} \Big(t_1,\x^{(1)}_{N_1}; \dots;t_M,\x^{(M)}_{N_M} \Big) 
=\det_{\substack
{1 \leq j \leq N_{m}, 1 \leq k \leq N_{n}, \\
1 \leq m, n \leq M}
}
\Bigg[
\mbK^{\xi}(t_m, x_{j}^{(m)}; t_n, x_{k}^{(n)} )
\Bigg],
$$
$0 \leq t_1 < \cdots < t_M < t_{*}$,
$1 \leq m \leq M$, 
$1 \leq N_m \leq N$,
$\x^{(m)}_{N_m} \in \A^{2 \pi r}_{N_m}, 1 \leq m \leq M \in \N$.
Here the integral kernel $\mbK^{\xi}: ([0, t_{*}) \times [0, 2 \pi r))^2 \mapsto \R$ is called
the spatio-temporal correlation kernel.
\end{df}

By Theorem 1.3 in \cite{Kat14},
DMR given by Proposition \ref{thm:DMR}
leads to the following result.

\begin{thm}
\label{thm:kernel}
For $\xi=\sum_{j=1}^N \delta_{u_j} \in \mM_0$ with $\u \in \check{\A}_N^{2 \pi r}$, 
the process $(\Xi(t), t \in [0, t_{*}), \P^{\xi})$ is
determinantal for any observable with the correlation kernel
\begin{eqnarray}
\mbK^{\xi}(s, x; t, y)
&=& \int_0^{2 \pi r} \xi(du) \,
p^r_{\sp(N)}(s,x|u) \cM_{\xi}^{u}(t,y)-\1(s>t) p^r_{\sp(N)}(s-t,x|y),
\nonumber\\
&& \hskip 3cm
(s,x), (t,y) \in [0, t_{*}) \times [0, 2 \pi r),
\label{eqn:K1}
\end{eqnarray}
where $p^r_{\sp(N)}$ and $\cM_{\xi}^u$ were
defined by (\ref{eqn:pkappa}) and (\ref{eqn:cM1}),
respectively.
\end{thm}
\vskip 0.3cm
\noindent{\bf Remark 4.} 
In the previous paper \cite{Kat15}, we showed that
eDYS$^{(2)}$ is determinantal only when it starts from
the configuration $\v$ with equidistant spacing in $[0, 2 \pi r)$.
In the present paper, we have given a new construction of
eDYS$^{(2)}$ and generalized the result.
\vskip 0.3cm

\SSC
{Concluding Remarks \label{sec:concluding_remarks}}

We now discuss a conjectural extension of
Theorem \ref{thm:kernel} to arbitrary initial configuration $\xi$
which can have multiple points in general. 
For $\xi = \sum_{j=1}^N \delta_{u_j}$, 
$s \in [0, t_*)$, $u, x \in [0, 2 \pi r)$, $z \in \C, \zeta \in \C \setminus \{u_{j}\}_{j=1}^N$,
$N \in \{2,3, \dots\}$, let
\begin{eqnarray}
\phi_{\xi}^u((s, x); z, \zeta)
&=& \frac{p^r_{\sp(N)}(s, x|\zeta)}{p^r_{\sp(N)}(s, x|u)} \times
\frac{\vartheta_1^{\prime}(0; N \tau(0))}{2 \pi r \vartheta_1((z-\zeta)/2 \pi r; N \tau(0))}
\nonumber\\
&\times& 
\frac{\vartheta_1((\overline{u}+z-u)/2\pi r; N \tau(0))}{\vartheta_1(\overline{u}/2 \pi r; N \tau(0))}
\prod_{\ell=1}^N \frac{\vartheta_1((z-u_{\ell})/2 \pi r; N \tau(0))}{\vartheta_1((\zeta-u_{\ell})/2 \pi r; N \tau(0))}.
\label{eqn:phiA1}
\end{eqnarray}
Note that
$$
\vartheta_1 \left( \frac{\zeta-u_{\ell}}{2 \pi r}; N \tau(0) \right)
=\frac{\zeta-u_{\ell}}{2 \pi r} \vartheta_1^{\prime}(0; N \tau(0))
+ {\cal O}((\zeta-u_{\ell})^2),
$$
in the vicinity of $\zeta=u_{\ell}, 1 \leq \ell \leq N$,
if $\xi \in \mM_0$. Thus, as a function of $\zeta$,
(\ref{eqn:phiA1}) has simple poles at $\zeta=u_{\ell}, 1 \leq \ell \leq N$,
when $\xi \in \mM_0$. 
Note that
$\vartheta_1^{\prime}(0; N \tau(0))
= 2 \pi e^{-N t_*/8 r^2} \prod_{n=1}^{\infty} (1-e^{-n N t_*/r^2})^3$
by (\ref{eqn:theta_prime1}).
Then we define
$$
\Phi_{\xi}^u((s,x); z)
= \frac{1}{2 \pi i} \oint_{\rC(\delta_u)} d \zeta \,
\phi_{\xi}^u((s,x); z, \zeta),
$$
where $\rC(\delta_u)$ is a closed contour on the complex plane $\C$
encircling the point $u$ once in the positive direction.
This is an entire function of $z$, and thus
\begin{eqnarray}
\cM_{\xi}^u((s, x)|(t, y))
&\equiv& \widetilde{\rE} [ \Phi_{\xi}^u((s,x); y+i \widetilde{b}(t)) ]
\nonumber\\
&=& \int_{\R} d \widetilde{b} \, 
\frac{e^{-\widetilde{b}^2/2t}}{\sqrt{2 \pi t}}
\Phi_{\xi}^u((s,x); y+i \widetilde{b}) 
\nonumber
\end{eqnarray}
gives a martingale, if we put $y=Y_{\sp(N)}(t), t \in [0, 2 \pi t_*)$
for any $u \in \supp \xi$.

For $\xi \in \mM$ with $\xi([0, 2\pi r))=N < \infty$, let
$\xi_*= \sum_{u \in \supp \xi} \delta_u$.
If $\xi \in \mM_0$, $\xi_*=\xi$ and it is easy to verify that
(\ref{eqn:K1}) is written as
\begin{eqnarray}
\mbK^{\xi}(s, x; t, y)
&=& \int_0^{2 \pi r} \xi_*(du) \,
p^r_{\sp(N)}(s,x|u) \cM_{\xi}^{u}((s,x) | (t, y)) -\1(s>t) p^r_{\sp(N)}(s-t,x|y),
\nonumber\\
&& \hskip 3cm
(s,x), (t,y) \in [0, t_{*}) \times [0, 2 \pi r).
\label{eqn:K2}
\end{eqnarray}
This function is well-defined as a spatio-temporal kernel
for general $\xi \in \mM$.
As an extension of Proposition 2.1 given in \cite{KT10}
for the original Dyson model with $\beta=2$, we expect the following.

\begin{conj}
\label{thm:conjecture1}
For any initial configuration $\xi \in \mM$ with $\xi([0, 2 \pi r)) = N < \infty$
and $\int_0^{2 \pi r} \xi(du) u - \kappa_N \in (0, 2 \pi r)$, 
the process $(\Xi(t), t \in [0, t_{*}), \P^{\xi})$ is
determinantal for any observable with the correlation kernel (\ref{eqn:K2}).
\end{conj}

It is an interesting future problem to construct
the infinite-particle systems of eDYS$^{(\beta)}$, $\beta >0$
by taking proper scaling limits $N \to \infty, r \to \infty$.
When $\beta=2$ the process is determinantal and
this problem will be reduced to the study concerning
the scaling limits of the spatio-temporal correlation kernel
given by (\ref{eqn:K1}) and (\ref{eqn:K2}).

\vskip 1cm
\noindent{\bf Acknowledgements} \quad
The present author would like to thank K. Takemura for useful comments
on the Calogero-Sutherland-Moser models.
A part of this work was done during the workshop
``New approaches to non-equilibrium and random systems: KPZ integrability, universality, applications and experiments"
(Jan.11--Mar.11, 2016) 
organized by Ivan Corwin, Pierre Le Doussal, and Tomohiro Sasamoto
at Kavli Institute for Theoretical Physics, University of California, Santa Barbara.
This research was supported in part by the National Science Foundation under Grant No.NSF PHY11-25915.
This work was also supported in part by
the Grant-in-Aid for Scientific Research
(C) (No.26400405),
(B) (No.26287019), and
(S) (No.16H06338)
 of Japan Society for the Promotion of Science.

\vskip 2cm
\appendix
\begin{LARGE}
{\bf Appendices}
\end{LARGE}

\SSC
{Notations and Formulas of Jacobi's Theta Functions
and Related Functions \label{sec:appendixA}}
Let
$$
z=e^{v \pi i}, \quad q=e^{\tau \pi i},
$$
where $v \in \C$ and $\Im \tau > 0$. 
The Jacobi theta functions are defined by follows \cite{WW27,NIST10}, 
\begin{eqnarray}
\vartheta_0(v; \tau) &=& 
\sum_{n \in \Z} (-1)^n q^{n^2} z^{2n} \nonumber\\
&=& 1+ 2 \sum_{n=1}^{\infty}(-1)^n e^{\tau \pi i n^2} \cos(2 n \pi v),
\nonumber\\
\vartheta_1(v; \tau) &=& i \sum_{n \in \Z} (-1)^n q^{(n-1/2)^2} z^{2n-1}
\nonumber\\
&=& 2 \sum_{n=1}^{\infty} (-1)^{n-1} e^{\tau \pi i (n-1/2)^2} \sin\{(2n-1) \pi v\},
\nonumber\\
\vartheta_2(v; \tau) 
&=& \sum_{n \in \Z} q^{(n-1/2)^2} z^{2n-1}
\nonumber\\
&=& 2 \sum_{n=1}^{\infty} e^{\tau \pi i (n-1/2)^2} \cos \{(2n-1) \pi v\},
\nonumber\\
\vartheta_3(v; \tau) 
&=& \sum_{n \in \Z} q^{n^2} z^{2n}
\nonumber\\
&=& 1 + 2 \sum_{n=1}^{\infty} e^{\tau \pi i n^2} \cos (2 n \pi v).
\label{eqn:theta}
\end{eqnarray}
(Note that the present functions 
$\vartheta_{\mu}(v; \tau), \mu=1,2,3$ are denoted by
$\vartheta_{\mu}(\pi v,q)$,
and $\vartheta_0(v;\tau)$ by $\vartheta_4(\pi v,q)$ in \cite{WW27}.)
With respect to $v$, 
$\vartheta_1(v; \tau)$ is odd and
$\vartheta_{\mu}(v; \tau)$, $\mu=0,2,3$ are even.
We see
\begin{eqnarray}
&& \vartheta_0(v; \tau) \sim 1, \quad
\vartheta_1(v; \tau) \sim 2 e^{\tau \pi i/4} \sin (\pi v), \quad
\vartheta_2(v; \tau) \sim 2 e^{\tau \pi i/4} \cos(\pi v), \quad
\vartheta_3(v; \tau) \sim 1,
\nonumber\\
&&\qquad \qquad \qquad \mbox{as} \quad
\Im \tau \to + \infty \quad
({\rm i.e.} \quad q=e^{\tau \pi i} \to 0).
\label{eqn:theta_asym}
\end{eqnarray}

For $\Im \tau >0$, $\vartheta_{\mu}(v; \tau)$, $\mu=0,1,2,3$
are holomorphic for $|v| < \infty$
and satisfy the partial differential equation
\begin{equation}
\frac{\partial \vartheta_{\mu}(v; \tau)}{\partial \tau}=
\frac{1}{4 \pi i} \frac{\partial^2 \vartheta_{\mu}(v; \tau)}{\partial v^2}.
\label{eqn:Jacobi_eq}
\end{equation}
The functions $\vartheta_{\mu}(v; \tau)$, $\mu=0,2,3$ are
expressed by $\vartheta_1(v; \tau)$ as
\begin{eqnarray}
\vartheta_0(v; \tau) &=& -i e^{\pi i (v+\tau/4)} \vartheta_1 \left( v+\frac{\tau}{2}; \tau \right),
\nonumber\\
\vartheta_2(v; \tau) &=& \vartheta_1 \left( v+\frac{1}{2}; \tau \right),
\nonumber\\
\vartheta_3(v; \tau) &=& e^{\pi i (v+\tau/4)}
\vartheta_1 \left( v+\frac{1+\tau}{2}; \tau \right).
\label{eqn:theta_relations}
\end{eqnarray}
They have the quasi-periodicity: for instance, $\vartheta_1$ satisfies
\begin{eqnarray}
\vartheta_1(v+1; \tau)
&=& -\vartheta_1(v; \tau),
\nonumber\\
\vartheta_1(v+\tau; \tau)
&=&-\frac{1}{z^2 q} \vartheta_1(v; \tau)
=-e^{-\pi i (2v+\tau)} \vartheta_1(v; \tau).
\label{eqn:quasi_periodic}
\end{eqnarray}
We can show that
\begin{equation}
{\vartheta_1}^{\prime}(0; \tau) \equiv
\left. \frac{d \vartheta_1(v; \tau)}{d v} \right|_{v=0}
= 2 \pi q^{1/4} \prod_{n=1}^{\infty} (1-q^{2n})^3.
\label{eqn:theta_prime1}
\end{equation}

The following functional equalities are known as
Jacobi's imaginary transformations \cite{WW27,NIST10},
\begin{eqnarray}
\vartheta_0(v; \tau)
&=& e^{\pi i/4} \tau^{-1/2} e^{-\pi i v^2/\tau}
\vartheta_2 \left( \frac{v}{\tau}; - \frac{1}{\tau} \right),
\nonumber\\
\vartheta_1(v; \tau)
&=& e^{3 \pi i/4} \tau^{-1/2} e^{-\pi i v^2/\tau}
\vartheta_1 \left( \frac{v}{\tau}; - \frac{1}{\tau} \right),
\nonumber\\
\vartheta_3(v; \tau)
&=& e^{\pi i/4} \tau^{-1/2} e^{-\pi i v^2/\tau}
\vartheta_3 \left( \frac{v}{\tau}; - \frac{1}{\tau} \right). 
\label{eqn:Jacobi_imaginary}
\end{eqnarray}

The Weierstrass elliptic functions
$\zeta$ and $\wp$ are
defined by
\begin{eqnarray}
\zeta(z) &=& \zeta(z| 2 \omega_1, 2 \omega_3)
\nonumber\\
\label{eqn:zeta1}
&=& \frac{1}{z}
+\sum_{(m,n) \in \Z^2 \setminus \{(0,0)\}}
\left[ \frac{1}{z-\Omega_{m,n}}+\frac{1}{\Omega_{m,n}}
+\frac{z}{{\Omega_{m,n}}^2} \right], \\
\wp(z) &=& \wp(z| 2 \omega_1, 2 \omega_3)
\nonumber\\
\label{eqn:wp1}
&=& \frac{1}{z^2}
+\sum_{(m,n) \in \Z^2 \setminus \{(0,0)\}}
\left[ \frac{1}{(z-\Omega_{m,n})^2}-\frac{1}{{\Omega_{m,n}}^2} \right],
\end{eqnarray}
where $\omega_1$ and $\omega_3$ are fundamental periods,
$\Omega_{m,n} = 2 m \omega_1+2n \omega_3$, and
\begin{eqnarray}
\eta_1(2 \omega_1, 2 \omega_3)
&=& \zeta(\omega_1 | 2 \omega_1, 2 \omega_3) 
\nonumber\\
\label{eqn:eta1}
&=& \frac{\pi^2}{\omega_1}
\left( \frac{1}{12} - 2 \sum_{n=1}^{\infty}
\frac{n q^{2n}}{1-q^{2n}} \right),
\end{eqnarray}
with $q=e^{\tau \pi i}$, $\tau=\omega_3/\omega_1$.
Note that $\zeta(z)$ is odd and $\wp(z)$ is even
as functions of $z$.
The following equality holds
(see Lemma 2.1 in \cite{Kat15}).
For $a, b, c \in \C$, 
\begin{eqnarray}
&& \zeta(a-b)\zeta(a-c)+\zeta(b-a)\zeta(b-c)+\zeta(c-a)\zeta(c-b)
\nonumber\\
&& \quad = \frac{1}{2} \Big\{
\zeta(a-b)^2+\zeta(b-c)^2+\zeta(a-c)^2 \Big\}
-\frac{1}{2} \Big\{ \wp(a-b)+\wp(b-c)+\wp(a-c) \Big\}.
\label{eqn:zeta_wp1}
\end{eqnarray}

The Dedekind modular function $\eta(\tau)$ is defined by
(see, for instance, Sec. 23.15 of \cite{NIST10})
\begin{equation}
\eta(\tau)=e^{\tau \pi i/12}
\prod_{n=1}^{\infty} (1-e^{2 n \tau \pi i}),
\quad \Im \tau > 0.
\label{eqn:Dedekind1}
\end{equation}
Its logarithmic derivative is related to $\eta_1$ defined by (\ref{eqn:eta1}) as
\begin{equation}
\frac{d}{d \tau} \log \eta(\tau)
= i \frac{\omega_1}{\pi}
\eta_1(2 \omega_1, 2 \tau \omega_1).
\label{eqn:Dedekind2}
\end{equation}


\SSC
{Basic Properties of $U_N(t, x)$ \label{sec:appendixB}}

The function $U_N(t, x), t \in [0, t_*), x \in (0, 2 \pi r)$
is defined by (\ref{eqn:cW1}).
The spatial derivatives have the following expressions 
\cite{WW27,NIST10},
\begin{eqnarray}
U_N^{\prime}(t,x)
&\equiv& \frac{\partial U_N}{\partial x}(t, x) 
= - \frac{1}{2 \pi r} \frac{\vartheta'(x/2 \pi r; N \tau(t))}
{\vartheta_1(x/2 \pi r; N \tau(t)) }
\nonumber\\
&=& -\zeta ( x | 2 \pi r,  2 N \omega_3(t) )
+\frac{x}{\pi r} \eta_1 ( 2 \pi r,  2 N \omega_3(t) )
\nonumber\\
\label{eqn:Wdx1}
&=& - \frac{1}{2 r} \cot \left( \frac{x}{2r} \right)
- \frac{2}{r} \sum_{n=1}^{\infty}
\frac{e^{-n N(t_*-t)/r^2}}{1-e^{-nN(t_*-t)/r^2}}
\sin \left( \frac{nx}{r} \right), \\
U_N^{\prime \prime}(t,x)
&\equiv& \frac{\partial^2 U_N}{\partial x^2}(t, x) 
=- \frac{1}{4 \pi^2 r^2}
\left\{ \frac{\vartheta_1^{\prime \prime}(x/2 \pi r^2; N \tau(t)) }
{\vartheta_1(x/2 \pi r; N \tau(t)) }
-\left( \frac{\vartheta_1^{\prime}(x/2 \pi r; N \tau(t)) }
{\vartheta_1(x/2 \pi r; N \tau(t))} \right)^2 \right\}
\nonumber\\
\label{eqn:Wd2x1}
&=& \wp ( x | 2 \pi r, 2 N \omega_3(t) )
+\frac{1}{\pi r} \eta_1 ( 2 \pi r, 2 N \omega_3(t) ),
\quad t \in [0, t_*), x \in (0, 2 \pi r). 
\end{eqnarray}
The above expressions imply the following symmetries,
\begin{equation}
U_N^{\prime}(t, -x) = - U_N^{\prime}(t, x),
\quad
U_N^{\prime \prime}(t, -x)=U_N^{\prime \prime}(t, x).
\label{eqn:BB}
\end{equation}
By (\ref{eqn:Jacobi_eq}), the temporal derivative of $U_N(t, x)$ is
given by
\begin{equation} 
\dot{U}_N(t,x) \equiv \frac{\partial U_N}{\partial t}(t, x) 
= \frac{i N}{2 \pi r^2}
\frac{\dot{\vartheta}_1(x/2 \pi r; N \tau(t))}
{\vartheta_1(x/2 \pi r; N \tau(t))}
= \frac{N}{8 \pi^2 r^2}
\frac{\vartheta_1^{\prime \prime}(x/2 \pi r; N \tau(t))}
{\vartheta_1(x/2 \pi r; N \tau(t))}.
\label{eqn:Wdt1}
\end{equation}
Then we obtain the following partial differential equation,
\begin{equation}
\dot{U}_N(t,x) = \frac{N}{2}
\Big( U_N^{\prime}(t,x)^2 - U_N^{\prime \prime}(t,x) \Big),
\quad t \in [0, t_*), \quad x \in (0, 2 \pi r).
\label{eqn:WNrel}
\end{equation}

\vskip 1cm



\end{document}